\PassOptionsToPackage{table}{xcolor}

\documentclass{article}

\usepackage{microtype}
\usepackage{graphicx}
\usepackage{subfigure}
\usepackage{booktabs} %
\usepackage{enumitem}
\usepackage{hyperref}

\newcommand{\comment}[1]{}

\usepackage[accepted]{icml2023}

\usepackage{amsmath}
\usepackage{amssymb}
\usepackage{mathtools}
\usepackage{amsthm}
\usepackage{bm}

\usepackage{thmtools,thm-restate}
\usepackage{mdframed}

\usepackage[capitalize,noabbrev]{cleveref}

\usepackage{latexsym}
\usepackage{calc}
\usepackage{listings}
\usepackage{mathrsfs} 
\usepackage{float}

\usepackage{caption}
\usepackage{multirow, hhline}

\usepackage{color}
\usepackage{pifont}
\definecolor{mydarkgreen}{RGB}{39,130,67}

\definecolor{mydarkred}{RGB}{192,47,25}

\newcommand{\cmark}{{\normalsize {\color{mydarkgreen}\ding{51}}}}%
\newcommand{\xmark}{{\normalsize {\color{mydarkred} \ding{55}}}}%

\theoremstyle{plain}
\newtheorem{theorem}{Theorem}[section]
\newtheorem{proposition}[theorem]{Proposition}
\newtheorem{lemma}[theorem]{Lemma}
\newtheorem{corollary}[theorem]{Corollary}
\newtheorem{example}[theorem]{Example}
\newmdtheoremenv{framedtheorem}[theorem]{Theorem}
\newmdtheoremenv{framedlemma}[theorem]{Lemma}
\newmdtheoremenv{framedexample}[theorem]{Example}
\newmdtheoremenv{framedassumption}[theorem]{Assumption}
\newmdtheoremenv{framedproposition}[theorem]{Proposition}

\theoremstyle{definition}
\newtheorem{definition}[theorem]{Definition}
\newtheorem{assumption}[theorem]{Assumption}
\theoremstyle{remark}
\newtheorem{remark}[theorem]{Remark}

\usepackage[textsize=tiny]{todonotes}

\newcommand{\mat}[1]{\mathbf{#1}}

\def\argmin{\mathop{\rm argmin}}
\def\Argmin{\mathop{\rm Argmin}}

\def\Def{\stackrel{\mathrm{def}}{=}}

\def\dom{{\rm dom \,}}
\def\beq{\begin{equation}}
\def\eeq{\end{equation}}

\newcommand{\nz}{{\rm nz\,}}

\def\R{\mathbb{R}}

\providecommand{\0}{\mathbf{0}}

\renewcommand{\aa}{\mathbf{a}}
\providecommand{\bb}{\mathbf{b}}

\providecommand{\ee}{\mathbf{e}}

\renewcommand{\gg}{\mathbf{g}}
\providecommand{\hh}{\mathbf{h}}

\renewcommand{\ss}{\mathbf{s}}
\providecommand{\tt}{\mathbf{t}}

\providecommand{\xx}{\mathbf{x}}
\providecommand{\yy}{\mathbf{y}}
\providecommand{\zz}{\mathbf{z}}

\providecommand{\cO}{\mathcal{O}}

\def\Mymod{{\rm \, mod \,}}

\def\BI{\begin{itemize}}
	\def\EI{\end{itemize}}

\renewcommand\arraystretch{1}
\def\ba{\begin{array}}
	\def\ea{\end{array}}
\def\beann{\begin{eqnarray*}}
	\def\eeann{\end{eqnarray*}}
\def\bea{\begin{eqnarray}}
\def\eea{\end{eqnarray}}

\def\BT{\begin{theorem}}
	\def\ET{\end{theorem}}
\def\BL{\begin{lemma}}
	\def\EL{\end{lemma}}
\def\BC{\begin{corollary}}
	\def\EC{\end{corollary}}
\def\BE{\begin{example}}
	\def\EE{\end{example}}
\def\BD{\begin{definition}}
	\def\ED{\end{definition}}
\def\BR{\begin{remark}}
	\def\ER{\end{remark}}
\def\BAS{\begin{assumption}}
	\def\EAS{\end{assumption}}
\def\BI{\begin{itemize}}
	\def\EI{\end{itemize}}
\def\BP{\begin{proposition}}
	\def\EP{\end{proposition}}

\def\BMP{\begin{minipage}{9.5cm}}
	\def\EMP{\end{minipage}}
\def\MPT{\begin{minipage}{11.5cm}}
	\def\EPT{\end{minipage}}

\def\la{\langle}
\def\ra{\rangle}

\def\QF{\hspace{5ex} \Box}

\icmltitlerunning{Second-order optimization with lazy Hessians}

\begin{document}

\twocolumn[
\icmltitle{Second-Order Optimization with Lazy Hessians}

\icmlsetsymbol{equal}{*}

\begin{icmlauthorlist}
\icmlauthor{Nikita Doikov}{equal,yyy}
\icmlauthor{El Mahdi Chayti}{equal,yyy}
\icmlauthor{Martin Jaggi}{yyy}
\end{icmlauthorlist}

\icmlaffiliation{yyy}{Machine Learning and Optimization Laboratory, EPFL, Switzerland}

\icmlcorrespondingauthor{Nikita Doikov}{nikita.doikov@epfl.ch}
\icmlcorrespondingauthor{El Mahdi Chayti}{el-mahdi.chayti@epfl.ch}
\icmlcorrespondingauthor{Martin Jaggi}{martin.jaggi@epfl.ch}

\icmlkeywords{Newton method, convex optimization, non-convex optimization, second-order methods}

\vskip 0.3in
]

\printAffiliationsAndNotice{\icmlEqualContribution} %

\begin{abstract}
We analyze Newton's method with lazy Hessian updates
for solving general possibly non-convex optimization problems.
We propose to reuse a previously seen Hessian for several iterations
while computing new gradients at each step of the method. 
This significantly reduces the overall arithmetic complexity
of second-order optimization schemes.
By using the cubic regularization technique, we establish fast global 
convergence of our method to a second-order stationary point, while the Hessian does not need to be
updated each iteration. 
For convex problems, we justify global and local superlinear rates
for lazy Newton steps with quadratic regularization, which is easier to compute.
The optimal frequency for updating the Hessian is once every~$d$ iterations, 
where $d$ is the dimension of the problem. This provably improves the total arithmetic 
complexity of second-order algorithms by a factor $\sqrt{d}$.
\end{abstract}

\section{Introduction}

\paragraph{Motivation.}

Second-order optimization algorithms are being widely used
for solving difficult ill-conditioned problems.
The classical Newton method approximates the objective function
by its second-order Taylor approximation at the current point.
A minimizer of this model serves as one step of the algorithm. Locally, Newton's method achieves very 
fast quadratic rate \cite{kantorovich1948functional},
when the iterates are in a neighborhood of the solution.
However, it may not converge if the initial point is far away
from the optimum.
For the global convergence of Newton's method, we need 
to use some of the various regularization techniques,
that have been intensively developed during recent decades
(including \textit{damped Newton steps} with line search \cite{kantorovich1948functional,ortega2000iterative,hanzely2022damped}, 
\textit{trust-region} methods \cite{conn2000trust}, 
\textit{self-concordant} analysis \cite{nesterov1994interior,bach2010self,sun2019generalized,dvurechensky2018global} or 
the notion of \textit{Hessian stability} \cite{karimireddy2018global},
\textit{cubic regularization}~\cite{nesterov2006cubic},
\textit{contracting-point} steps \cite{doikov2020convex,doikov2022affine}, 
or \textit{gradient regularization} techniques 
\cite{mishchenko2021regularized,doikov2021gradient}). In all of these approaches, the arithmetic complexity of each step remains costly, since it requires computing both the gradient and Hessian of the objective function,
and then perform some factorization of the matrix. 
Hessian computations are prohibitively expensive in any large-scale setting.

In this work, we explore a simple possibility for a significant acceleration of second-order schemes in terms of their total arithmetic complexity.
The idea is to keep a previously
computed Hessian (a stale version of the Hessian) for several iterations, while using fresh
gradients in each step of the method. 
Therefore, we can save a lot of computational resources
by reusing the old Hessians. We call this \textit{lazy Hessian updates}.
We justify that this idea can be successfully employed in most established second-order schemes, 
resulting in a provable improvement of the overall performance.

\begin{table}[h!]
	
	\centering
	\setlength{\arrayrulewidth}{0.5mm}
	\renewcommand{\arraystretch}{1.8}
	\setlength{\tabcolsep}{3pt}
	\arrayrulecolor[HTML]{66BB6A}
	
	\begin{tabular}{| c | c | c | c | c |}
		\hhline{-----}
		\cellcolor[HTML]{FFFDE7}
		$\boldsymbol{\nabla^2 f(\xx_0)}$ & 
		\multicolumn{3}{c|}{  
			\cellcolor[HTML]{FFFDE7}
			\large
			$\overset{\textbf{reuse Hessian}}{\boldsymbol{\longrightarrow}}$	 
		}  &
		\cellcolor[HTML]{FFFDE7}
		$\!\boldsymbol{\nabla^2 f(\xx_m)}\!$ \\
		\hhline{-----}
		\cellcolor[HTML]{FFFDE7}
		$\boldsymbol{ \nabla f(\xx_0) }$ &
		\cellcolor[HTML]{FFFDE7}
		$\boldsymbol{ \nabla f(\xx_1)} $ & 
		\cellcolor[HTML]{FFFDE7}
		$ \boldsymbol{ \ldots }$ & 
		\cellcolor[HTML]{FFFDE7}
		$\boldsymbol{\nabla f(\xx_{m - 1})}$ &
		\cellcolor[HTML]{FFFDE7}
		$\boldsymbol{\nabla f(\xx_{m})}$ 
		\\[-\arrayrulewidth]
		\hhline{-----}
	\end{tabular}
	\caption*{\emph{Lazy Hessian Updates:} compute a new Hessian once per $m$ iterations.}
	\centering
	
\end{table}

\begin{table*}[h!]
	\centering
	\small
	\renewcommand{\arraystretch}{1.6}
	\begin{tabular}{cccccc}
		\arrayrulecolor[HTML]{000000}
		\hline
		Method & 
		\hspace*{-12pt}
		\begin{tabular}{c}Global\\[-5pt] Complexity\end{tabular} & 
		\hspace*{-15pt}
		\begin{tabular}{c}Local\\[-5pt]superlinear\\[-5pt] convergence\end{tabular} & 
		\hspace*{-15pt}
		\begin{tabular}{c}Non-convex\\[-5pt] problems\end{tabular} &
		\hspace*{-15pt}
		\begin{tabular}{c}Hessian\\[-5pt] computations\end{tabular} & 
		 Reference \\
		\hline
		Gradient Method & 
		\hspace*{-12pt} $\cO(\varepsilon^{-2})$ &
		\hspace*{-15pt} \xmark & 
		\hspace*{-15pt} \cmark & 
		\hspace*{-15pt} \textbf{not used} & 
		\cite{ghadimi2013stochastic}  \\
		Classical Newton & 
		\hspace*{-15pt} \xmark & 
		\hspace*{-15pt} \cmark & 
		\hspace*{-15pt} \xmark & 
		\hspace*{-15pt} {\color{mydarkred}\textbf{every step}} & 
		\cite{nesterov2018lectures} \\
		Shamanskii Newton & 
		\hspace*{-15pt} \xmark & 
		\hspace*{-15pt} \cmark & 
		\hspace*{-15pt} \xmark & 
		\hspace*{-15pt} {\color{mydarkgreen}\textbf{once per $\boldsymbol{m}$ steps}}& 
		\cite{shamanskii1967modification} \\
		Cubic Newton & 
		\hspace*{-12pt} $\cO(\varepsilon^{-3/2})$ & 
		\hspace*{-15pt} \cmark & 
		\hspace*{-15pt} \cmark  & 
		\hspace*{-15pt} {\color{mydarkred}\textbf{every step}}& 
		\cite{nesterov2006cubic} \\
		\hline
		Cubic Newton with Lazy Hessians\;\; &  
		\hspace*{-12pt} $\cO(\sqrt{m} \varepsilon^{-3/2})$ & 
		\hspace*{-15pt} \cmark & 
		\hspace*{-15pt} \cmark & 
		\hspace*{-13pt} {\color{mydarkgreen}\textbf{once per $\boldsymbol{m}$ steps}}& 
		 \textbf{ours} (Algorithm \ref{alg:LazyCubicNewton}) \\
		Reg. Newton with Lazy Hessians\, &  
		\hspace*{-12pt} $\cO(\sqrt{m}  \varepsilon^{-3/2})$ & 
		\hspace*{-15pt} \cmark & 
		\hspace*{-15pt} \xmark & 
		\hspace*{-13pt} {\color{mydarkgreen}\textbf{once per $\boldsymbol{m}$ steps}}& 
		 \textbf{ours} (Algorithm \ref{alg:LazyGradRegNewton}) \\
		\hline
		
	\end{tabular}
	\caption{A comparison of our methods with the classical deterministic algorithms.
		The Global Complexity shows how much gradient computations are needed for a method to find 
		a point $\bar{\xx}$ with the small gradient norm: $\| \nabla f(\bar{\xx}) \|_{*} \leq \varepsilon$.}
	\label{TableComplexities}
\end{table*}

For general problems\footnote{%
	Assuming that the cost of computing one Hessian and its appropriate factorization
	is $d$ times the cost of one gradient step. See Examples~\ref{ex:separable}, \ref{ex:auto_diff},
	and \ref{ex:finite_diff}. }%
, the optimal schedule is $m := d$ (update the Hessian every $d$ iterations),
where $d$ is the dimension of the problem.
In these cases, we gain a provable computational acceleration of our methods by a factor $\sqrt{d}$. The key to our analysis is an increased value of the regularization constant. 
In a general case, this should be proportional to $mL$,
where $L$ is the Lipschitz constant of the Hessian. 
Thus, for the lazy Newton steps, we should regularize our model stronger to balance the possible errors
coming from the Hessian inexactness, and so the steps of the method become shorter.
Fortunately, we save much more in terms of overall computational cost, which makes our approach appealing.

\paragraph{Related Work.}

The idea of using an old Hessian for several steps of Newton's method
probably appeared for the first time in \cite{shamanskii1967modification},
with the study of local convergence of this process as applied to solving
systems of nonlinear equations. The author proved a local quadratic rate
for the iterates when updating the Hessian, and a linear rate with an improving factor for the iterates in between.
Later on, this idea has been successfully combined with
Levenberg-Marquardt \cite{levenberg1944method,marquardt1963algorithm}  
regularization approach in \cite{fan2013shamanskii},
with damped Newton steps in \cite{lampariello2001global,wang2006further},
and with the proximal Newton-type methods for convex optimization in \cite{adler2020new}.
In these papers, it was established that the methods possess 
an asymptotic global convergence to a solution, without explicit 
non-asymptotic bounds on the actual rate.  Compared with these works, our new second-order algorithms 
with the lazy Hessian updates use a different analysis
that is based on the modern globalization techniques
(\textit{cubic regularization} \cite{nesterov2006cubic} and 
\textit{gradient regularization} \cite{mishchenko2021regularized,doikov2021gradient}).
It allows equipping our methods with provably fast global rates for wide classes of convex and non-convex optimization problems.
Thus, we prove a global complexity $\cO(1 / \varepsilon^{3/2})$ of the lazy regularized Newton steps
for achieving a second-order stationary point with the gradient norm bounded by $\varepsilon$,
while saving a lot of computational resources by reusing the old Hessian 
(see Table~\ref{TableComplexities} for the comparison of the global complexities). Another interesting connection can be established to recently developed distributed Newton-type methods \cite{islamov2022distributed},
where the authors propose to use a probabilistic aggregation and compression of the Hessians, as for example in federated learning.
In the particular case of the single node setup,
their algorithm also evaluates the Hessian rarely which is similar in spirit to our approach
but using different aggregation strategies. 
In their technique, the node needs to consider the Hessian each iteration. Then, using a certain criterion,
it decides whether to employ the new Hessian for the next step or to keep using just the old information,
while in our approach the Hessian is computed only \textit{once} per $m$ iterations. 
The authors proved local linear and superlinear rates that are independent on the condition number.
At the same time, in our paper, we primarily focus
on global complexity guarantees and thus establish a provable computational improvement
by a factor $\sqrt{d}$ for our methods as compared to updating the Hessian in every step.
We also justify local superlinear rates for our approach.

\paragraph{Contributions.}
We develop new efficient second-order optimization algorithms
and equip them with the global complexity guarantees. More specifically,
\vspace{-2mm}
\begin{itemize}[noitemsep]
	\item We propose the lazy Newton step with cubic regularization 
	(Section~\ref{SectionLazyCubic}).
	 It uses the gradient computed at the current point and the Hessian at some different point from
	the past trajectory. 
	We quantify the error effect coming from the inexactness in the second-order information and
	formulate the progress for one method step (Theorem~\ref{ThProgress}).
	We show how to balance the errors for $m$ consecutive lazy Newton steps,
	by increasing the regularization parameter proportionally. 
	
	\item Based on that, we develop the Cubic Newton with Lazy Hessians  (Algorithm~\ref{alg:LazyCubicNewton}) 
	and establish its fast global convergence to a \textit{second-order stationary point}
	 (Theorem~\ref{ThGlobCompl} in Section~\ref{SectionGlobalCubic}).
	 This avoids that
	 our method could get stuck in saddle points.
	 In the case $m := 1$ (updating
	the Hessian each iteration), our rate recovers the classical rate of full Cubic Newton \cite{nesterov2006cubic}.
	Then,
	we show that taking into account the actual arithmetic cost of the Hessian computations, the
	\textit{optimal} choice for one phase of the method is $m := d$, which improves the total complexity of
	Cubic Newton by a factor $\sqrt{d}$ (see Corollary \ref{CorOptimal}).

	\item We show how to improve our method when the problem is \textit{convex} in Section~\ref{SectionConvex}.
	We develop the Regularized Newton with Lazy Hessians (Algorithm~\ref{alg:LazyGradRegNewton}),
	which replaces the cubic regularizer in the model by \textit{quadratic} one.  This makes the subproblem much easier to solve,
	involving just one standard matrix inversion, while keeping the fast rate of the original cubically
	regularized method. 
	
	\item We study the local convergence of our new algorithms
	in Section~\ref{SectionLocal} (see Theorem~\ref{ThLocal} and Theorem~\ref{ThGRLocal}). 
	We prove that they both enjoy superlinear convergence rates. 
	As a particular case, we also justify the local quadratic convergence for the \textit{classical Newton method} 
	(without regularization)
	with the lazy Hessian updates.
	
	\item Illustrative numerical experiments are provided.
\end{itemize}

\vspace{-1mm}
\paragraph{Notation.}

We consider the following optimization problem,
\beq \label{MainProblem}
\ba{c}
\min\limits_{\xx \in \R^d} f(\xx),
\ea
\eeq
where $f$ is a several times differentiable function, not necessarily convex, that we assume to be bounded from below: $\inf\limits_{\xx} f(\xx) \geq f^\star$.
We denote by $\nabla f(\xx) \in \R^d$ its gradient
and by $\nabla^2 f(\xx) \in \R^{d \times d}$ its Hessian,
computed at some $\xx \in \R^d$.

Our analysis will apply both if optimization is over the standard Euclidean space, and also more generally if a different norm is defined by an arbitrary fixed symmetric positive definite matrix $\mat{B} = \mat{B}^\top \succ 0$. In that case the norm becomes\vspace{-1mm}
\beq \label{EuclideanNorm}
\ba{rcl}
\|\xx\| & \Def & \la \mat{B}\xx, \xx \ra^{1/2}, \qquad \xx \in \R^d.
\ea
\eeq
Thus, the matrix $\mat{B}$ is responsible for fixing the coordinate system in our problem.
In the simplest case, we can choose $\mat{B} := \mat{I}$ (identity matrix),
which recovers the classical Euclidean geometry.
A more advanced choice of $\mat{B}$ can take into account a specific structure of the problem
(see Example~\ref{ex:separable}).
The norm for the dual objects (gradients) is defined in the standard way,
\vspace{-2mm}
\[
\ba{rcl}
\|\gg\|_{*} & \!\! \Def \!\! & \sup\limits_{\xx : \, \|\xx\| \leq 1} \la \gg, \xx \ra
\; = \;
\la \gg, \mat{B}^{-1} \gg \ra^{1/2}, \quad \gg \in \R^d.
\ea
\]
\vspace{-2mm}
The induced norm for a matrix $\mat{A} \in \R^{d \times d}$ is defined by
\[
\ba{rcl}
\|\mat{A}\| &  \Def \!\! & \sup\limits_{\substack{\xx, \yy: \\ \|\xx\| \leq 1, \|\yy\| \leq 1}}
\la \mat{A}\xx, \yy\ra
\; = \;
\sup\limits_{\xx: \, \|\xx\| \leq 1} \|\mat{A}\xx\|_{*}.
\ea
\]
We assume that the Hessian of $f$ is Lipschitz continuous, with some constant $L > 0$:
\vspace{-1mm}
\beq \label{OldLipHess}
\ba{rcl}
\!\!\!
\| \nabla^2 f(\xx) - \nabla^2 f(\yy) \| & \!\!\! \leq \!\!\! & L \|\xx - \yy\|,
\;\; \forall \xx, \yy \in \R^d.
\ea
\eeq

\newpage
\section{Lazy Newton Steps}
\label{SectionLazyCubic}

Let us introduce the following lazy Newton steps,
where we allow the gradient and Hessian to be computed at different points $\xx, \zz \in \R^d$.
Thus, we denote the following quadratic model of the objective:
$$
\ba{c}
\!\! Q_{\xx, \zz}(\yy)  \Def 
\la \nabla f(\xx), \yy - \xx \ra 
\! + \!\frac{1}{2}\la \nabla^2 f(\zz)(\yy - \xx), \yy - \xx \ra.
\ea
$$
We use cubic regularization of our model, for some $M > 0$:
\beq \label{StepDef}
\ba{c}
\!\!\! \mat{T}_{M}(\xx, \zz)  \in  \Argmin\limits_{\yy \in \R^d} 
\Bigl\{ 
Q_{\xx, \zz}(\yy)
+ \frac{M}{6}\|\yy - \xx\|^3
\Bigr\}.
\ea
\eeq

For $\zz := \xx$ this is the 
iteration of Cubic Newton~\cite{nesterov2006cubic}. 
Thus, in our scheme, we can reuse the Hessian from a previous step $\zz$ without recomputing it,
which significantly reduces the overall iteration cost.  

Our definition implies that the point $\mat{T} = \mat{T}_{M}(\xx, \zz)$
is a \textit{global minimum} of the cubically regularized model,
which is generally non-convex. However, it turns out that we can compute this point efficiently
by using standard techniques developed initially for trust-region methods~\cite{conn2000trust}. 
Let us denote $r \Def \|\mat{T} - \xx\|$. 
The solution to the subproblem~\eqref{StepDef} satisfies the following stationarity conditions:
$$
\ba{c}
\begin{cases}
&\nabla f(\xx) + \nabla^2 f(\zz)(\mat{T} - \xx) + \frac{Mr}{2} \mat{B}(\mat{T} - \xx)
\;\; = \;\; \0, \\
& \nabla^2 f(\zz) + \frac{Mr}{2} \mat{B} \;\; \succeq \;\; 0.
\end{cases}
\ea
$$

Thus, in the non-degenerate case, one step can be represented in the following form: 
\beq \label{OldCubicStep}
\ba{rcl}
\mat{T} & = & \xx - \bigl( \nabla^2 f(\zz) + \frac{Mr}{2} \mat{B} \bigr)^{-1} \nabla f(\xx),
\ea
\eeq
and the value $r > 0$ can be found by solving the corresponding univariate nonlinear equation
\citep[Section 5]{nesterov2006cubic}.
It can be done very efficiently from a precomputed \textit{eigenvalue} or the \textit{tridiagonal}
decomposition of the Hessian. Typically, it is of similar cost as for matrix inversion in the classical Newton step. We discuss the computation of the iterate $\mat{T}_{M}(\xx, \zz)$ in more details in Section~\ref{SubsectionFactorization}.
Let us define the following quantity, for $\yy \in \R^d$:
\beq \label{OldXiDef}
\ba{rcl}
\xi(\yy) & \Def & \Bigl[  -\lambda_{\min}\bigl( 
\mat{B}^{-1/2} \nabla^2 f(\yy) \mat{B}^{-1/2}
\bigr) \Bigr]_{+},
\ea
\eeq
where $[t]_{+} \Def \max\{ t, 0 \}$ denotes the positive part,
and $\lambda_{\min}(\cdot)$ is the smallest eigenvalue of a symmetric matrix.
If $\nabla^2 f(\yy) \succeq 0$ for a certain $\yy \in \R^d$, then $\xi(\yy) = 0$.
Otherwise, $\xi(\yy)$ shows how big (in absolute value) 
the smallest eigenvalue of the Hessian is w.r.t. a fixed matrix $\mat{B} \succ 0$.

\newpage
\begin{framedtheorem} \label{ThProgress}
Let $M \geq L$. Then, for one cubic step~\eqref{StepDef}, it holds
\beq \label{OneStepProgress}
\ba{cl}
& \!\!\!\!\!\!\! f(\xx) - f(\mat{T}) \\
\\
& \!\!\!\!\!\!\! \geq \; 
\max \Bigl\{ \frac{1}{648M^2}\xi(\mat{T})^3, 
\frac{1}{72 \sqrt{2M}} \| \nabla f(\mat{T}) \|_{*}^{3 / 2}
\Bigr\} \\
\\
& \!\!\!\!\!\!\! \;\;\;\; + \; \frac{M}{48} r^3
- \frac{11L^3}{M^2}\|\zz - \xx\|^3.
\ea
\eeq
\end{framedtheorem}

Theorem~\ref{ThProgress} shows how much progress we can expect from
one lazy Newton step with cubic regularization.
The price for using the lazy Hessian is the last term in the progress bound \eqref{OneStepProgress},
which vanishes if $\zz := \xx$ (updating the Hessian for the current iterate).

\section{Global Convergence Rates}
\label{SectionGlobalCubic}

Let us consider one phase of the potential algorithm when we compute the Hessian once at 
a current point $\zz := \xx_0$ and then perform $m$ cubic steps~\eqref{StepDef} in a row:
\beq \label{mIters}
\xx_{i}  =  \mat{T}_M(\xx_{i-1}, \zz), 
\quad i=1,\dots,m.
\eeq
We use some fixed regularization constant $M > 0$.
From \eqref{OneStepProgress}, we see that the error from using the old Hessian 
is proportional to the cube of the distance to the starting point: $\| \xx_{i - 1} - \xx_{0} \|^3$,
for each $1 \leq i \leq m$.
Hopefully, we can balance off the \textit{accumulating error} by aggregating as well
the positive terms from \eqref{OneStepProgress},
and by choosing a sufficiently big value of the regularization parameter $M$.
Hence, our aim would be to ensure:
\beq \label{IneqToProve}
\ba{rcl}
\!\!
\frac{M}{48}  \sum\limits_{i = 1}^m \|\xx_i - \xx_{i - 1}\|^3
& \! \geq \! &
\frac{11 L^3}{M^2} \sum\limits_{i = 1}^m \| \xx_{i - 1} - \xx_{0} \|^3.
\ea
\eeq
In fact, we show (see Appendix~\ref{AppendixSectionGlobal}) that it is enough to choose $\boxed{M\geq 6 mL}$
to guarantee \eqref{IneqToProve}. Thus 
we deal with the accumulating error for $m$ consecutive lazy steps in \eqref{mIters}.

After each inner phase of cubic steps \eqref{mIters} with stale Hessian, we recompute the Hessian at a new snapshot point $\zz$ in a new outer round. The length of each inner phase $m \geq 1$ is the key parameter
of the algorithm.
For convenient notation, we denote by $\pi(k)$ the highest multiple of $m$ which is less than or equal to $k$:
\beq \label{PiDef}
\ba{rcl}
\pi(k) & \Def & k - k \Mymod m, \qquad k \geq 0.
\ea
\eeq

The main algorithm is defined as follows:

\begin{algorithm}[h!]
	\caption{Cubic Newton with Lazy Hessians}
	\label{alg:LazyCubicNewton}
	\begin{algorithmic}
		\STATE {\bfseries Input:} $\xx_0 \in \R^d$, $m \geq 1$, $L > 0$. Choose $M \geq 0$.
		\FOR{$k=0,1,\dotsc$}
		\STATE Set last snapshot point $\zz_k = \xx_{\pi(k)}$
		\STATE Compute lazy cubic step $\xx_{k+1} = \mat{T}_{M}(\xx_k, \zz_k) $
		\ENDFOR
	\end{algorithmic}
\end{algorithm}

As suggested by our previous analysis, we will use a simple fixed rule for the regularization parameter:
\beq \label{MCubicNewton}
\boxed{
	\ba{rcl}
	M & := & 6mL
	\ea
}
\eeq

\begin{framedtheorem} \label{ThGlobCompl}
Let $M$ be fixed as in \eqref{MCubicNewton}. Assume that the gradients
for previous iterations $\{ \xx_i \}_{i = 1}^k$
of Algorithm~\ref{alg:LazyCubicNewton}  are higher than a desired error level $\varepsilon>0$:\vspace{-2mm}
\beq \label{OldGradBounded}
\ba{rcl}
\| \nabla f(\xx_i) \|_{*} & \geq & \varepsilon.
\ea
\eeq
Then, the number of iterations to reach  accuracy $ \| \nabla f(\xx_{k + 1}) \|_{*} \leq \varepsilon$  is at most
\beq \label{OldBoundAllIters}
\ba{rcl}
k & \leq & \cO\Bigl( \frac{\sqrt{m L} (f(\xx_0) - f^\star) }{\varepsilon^{3/2}} \Bigr),
\ea
\eeq
The total number of Hessian updates  during these iterations is 
\vspace{-4mm} 
\beq \label{BoundAllHess}
\ba{rcl}
t & \leq & \cO\Bigl( \frac{\sqrt{L} (f(\xx_0) - f^\star)}{\varepsilon^{3/2} \sqrt{m}} \Bigr).
\ea
\eeq
For the minimal eigenvalues of all Hessians, it holds that
\vspace{-4mm}
\beq \label{HessEigenBound}
\ba{rcl}
 \min\limits_{1 \leq i \leq k}
\xi(\xx_i) 
& \leq &
 \cO\Bigl(\frac{M^2(f(\xx_0) - f^\star)}{k} \Bigr)^{1/3}.
\ea
\eeq
\end{framedtheorem}

Let us assume that the desired accuracy level is small, i.e.
$$
\ba{rcl}
\varepsilon & \leq &  (f(\xx_0) - f^\star)^{2/3} \bigl(\frac{L}{m}\bigr)^{1/3},
\ea
$$
which requires from the method to use \textit{several Hessians}.

To choose the parameter $m$, we need to take into account the 
actual computational efforts required in all iterations of our method.
We denote by {\tt GradCost}  the arithmetic complexity 
of one gradient step, and by {\tt HessCost} the cost of computing the
Hessian and its appropriate factorization.
In general, we have
\beq \label{HessGradCost}
\ba{rcl}
\text{\tt HessCost} & = & d \cdot \text{\tt GradCost},
\ea
\eeq
where $d$ is the dimension of our problem \eqref{MainProblem}.

\BE
\label{ex:separable}
Let
$f(\xx) = \frac{1}{n} \sum\limits_{i = 1}^n \phi( \la \aa_i, \xx \ra )$,
where $\aa_i \in \R^d$, $1 \leq i \leq n$ are given data
and $\phi$ is a fixed univariate smooth function.
In this case,
$$
\ba{rcl}
\nabla f(\xx) & \! = \! & \mat{A}^\top \ss(\xx),
\;\;\; \text{and} \;\;\;\;
\nabla^2 f(\xx) \; = \; \mat{A}^\top \mat{Q}(\xx) \mat{A},
\ea
$$
where $\mat{A} \in \R^{n \times d}$ is the matrix with rows $\aa_1, \ldots, \aa_n$;
$\ss(\xx) \in \R^n$ and $\mat{Q}(\xx) \in \R^{n \times n}$ 
is a diagonal matrix 
given by
\vspace{-2mm}
$$
\ba{rcl}
\bigl[ \ss(\xx) \bigr]_{(i)} & \! \Def \! & \frac{1}{n}\phi'( \la \aa_i, \xx \ra ), \\
\bigl[ \mat{Q}(\xx) \bigr]_{(i, i)} & \!\Def \! & \frac{1}{n}\phi''( \la \aa_i, \xx \ra ).
\ea
$$
Because of the Hessian structure, it is convenient to use the matrix $\mat{B}:= \mat{A}^\top \mat{A}$ 
to define our primal norm \eqref{EuclideanNorm}.
Indeed, for this choice of the norm, we can show that the Lipschitz constant of the Hessian depends only on the loss function $\phi$
(i.e. it is an absolute numerical constant which does not depend on data). 
At the same time, when using the identity matrix $\mat{I}$, the corresponding Lipschitz constant depends on the size of the data matrix $\|\mat{A}\|$.
Thus, in the latter case, the value of the Lipschitz constant can be very big.

Let us assume that the cost of computing $\phi'(\cdot)$
and $\phi''(\cdot)$ is $\cO(1)$
which does not depend on the problem parameters.
We denote by $\nz(\mat{A})$ the number of nonzero elements of $\mat{A}$.
Then,
\vspace{-2mm}
$$
\ba{rcl}
\text{\tt HessCost} & = & \cO\bigl( d \cdot \nz(\mat{A}) +  d^3 \bigr),
\ea
$$
where the last cubic term comes from computing
a matrix factorization (see Section~\ref{SubsectionFactorization}),
and
$$
\ba{rcl}
\text{\tt GradCost} & = & 
\cO\bigl(  \nz(\mat{A}) + d^2 \bigr),
\ea
$$
where the second term is from using the factorization.
Thus, relation~\eqref{HessGradCost} is satisfied.
\EE

\BE
\label{ex:auto_diff}
Let the representation of our objective be given by the computational graph (e.g. a neural network).
In this case, we can compute the Hessian-vector product $\nabla^2 f(\xx) \hh$
for any $\xx, \hh \in \R^d$ at the same cost as its gradient $\nabla f(\xx)$
by using automatic differentiation technique \citep[Chapter 7.2]{nocedal2006numerical}.
Thus, in general we can compute the Hessian as $d$ Hessian-vector products,
\vspace{-2mm}
$$
\ba{rcl}
\nabla^2 f(\xx) & = & \bigl[ \, \nabla^2 f(\xx)\ee_1 \, \big| \, \ldots \, \big| \, \nabla^2 f(\xx)\ee_d  \, \bigr],
\ea
$$
where $\ee_1, \ldots, \ee_d$ are the standard basis vectors.
This satisfies \eqref{HessGradCost}.
\EE

\BE
\label{ex:finite_diff}
For any differentiable function, we can use 
an approximation of its Hessian by finite differences.
Namely, for a fixed parameter $\delta > 0$, we can form
$\mat{H}_{x, \delta} \in \R^{d \times d}$ as
\vspace{-1mm}
$$
\ba{rcl}
\bigl[ \mat{H}_{x, \delta} \bigr]^{(i, j)}
& \!\! := \!\! & 
\frac{1}{\delta} \bigl[ 
\nabla f(\xx + \delta \ee_i) - \nabla f(\xx)
\bigr]^{(j)},
\ea
$$
where $\ee_1, \ldots, \ee_d$ are the standard basis vectors,
and use the following symmetrization as an approximation of the true Hessian
(see \cite{nocedal2006numerical,cartis2012oracle,grapiglia2022cubic}
for more details):
\vspace{-1mm}
$$
\ba{rcl}
\frac{1}{2} \Bigl[
\mat{H}_{x, \delta} + \mat{H}_{x, \delta}^{\top}
\Bigr] & \approx & \nabla^2 f(\xx).
\ea
$$
Thus, forming the Hessian approximation requires $d + 1$ gradient computations,
which is consistent with~\eqref{HessGradCost}.
\EE

Hence, according to Theorem~\ref{ThGlobCompl}, the total arithmetic complexity
of Algorithm~\ref{alg:LazyCubicNewton} can be estimated as
\vspace{-1mm}
\beq \label{ArithCompl}
\ba{cl}
& {\tt Arithmetic \;\; Complexity} \\
\\
& \;\;\;\; = \;\;\;\;
k \times {\tt GradCost}  \;\; + \;\;
t \times {\tt HessCost} \\
\\
& 
\!\!\! \overset{\eqref{OldBoundAllIters},\eqref{BoundAllHess},\eqref{HessGradCost}}{\leq} 
\cO\Bigl( \bigl(  \sqrt{m} + \frac{d}{\sqrt{m}}  \bigr)  
\frac{\sqrt{L} (f(\xx_0) - f^\star)}{\varepsilon^{3/2}} \Bigr) \\
\\
& \qquad \qquad 
\times \;\;
{\tt GradCost}.
\ea
\eeq

\BC 
For $\boxed{m = 1}$ we update Hessian at each step. It corresponds to the full Cubic Newton \cite{nesterov2006cubic}, and 
Theorem~\ref{ThGlobCompl} recovers its global iteration complexity:
\vspace{-2mm}
$$
\ba{rcl}
k & = & t \;\; \leq \;\; \cO\Bigl(  \frac{\sqrt{L}(f(\xx_0) - f^\star)}{\varepsilon^{3/2}} \Bigr).
\ea
$$
\vspace{-1mm}
Consequently, the total arithmetic complexity of the method is bounded by
\vspace{-1mm}
\beq \label{FullCubicNewton}
\ba{c}
\cO\Bigl(  d \cdot \frac{\sqrt{L}(f(\xx_0) - f^\star)}{\varepsilon^{3/2}} \Bigr) \cdot \tt{GradCost}.
\ea
\eeq
\EC

\BC \label{CorOptimal}
For $\boxed{m = d}$ we obtain the \underline{optimal choice} for the length of one phase,
which minimizes the right hand side of \eqref{ArithCompl}.
The total arithmetic complexity of Algorithm~\ref{alg:LazyCubicNewton} becomes
$$
\ba{c}
\cO\Bigl(  \sqrt{d} \cdot \frac{\sqrt{L}(f(\xx_0) - f^\star)}{\varepsilon^{3/2}} \Bigr) \cdot \tt{GradCost}.
\ea
$$
This improves upon the full Cubic Newton by factor $\sqrt{d}$.
\EC

Now, let us look at the minimal eigenvalue of all the Hessians at points generated by our algorithm.

\BC
Let us fix some $\varepsilon > 0$ and perform
\beq \label{KChoice}
\ba{rcl}
k & = & \frac{\sqrt{mL}(f(\xx_0) - f^\star)}{\varepsilon^{3/2}}
\ea
\eeq
iterations of Algorithm~\ref{alg:LazyCubicNewton}. According to Theorem~\ref{ThGlobCompl},
we thus ensure
$
\min\limits_{1 \leq i \leq k} \| \nabla f(\xx_i) \|_{*}  \leq  \cO(\varepsilon).
$
At the same time,
$$
\ba{c}
 \min\limits_{1 \leq i \leq k} \xi(\xx_i)
\overset{\eqref{HessEigenBound}}{\leq}
\Bigl(\frac{648M^2(f(\xx_0) - f^\star)}{k} \Bigr)^{1/3} 
\!
\overset{\eqref{KChoice}}{=}
2^{5/3} 3^2 \sqrt{mL\varepsilon}.
\ea
$$

Therefore, the negative eigenvalues of the Hessians cannot be big. 
For $\cO(\varepsilon)$ level of gradient norm
we guarantee $\cO(\sqrt{mL\varepsilon})$ level
 for the smallest eigenvalue \eqref{OldXiDef}.
\EC

\section{Minimizing Convex Functions}
\label{SectionConvex}

In this section, we assume that the objective 
in our problem~\eqref{MainProblem} is \textit{convex}.
Thus, all the Hessians are positive semidefinite, i.e.,
$\nabla^2 f(\xx)  \succeq  0,\ $ $\forall \xx\in \R^d$.

Then, we can apply the \textit{gradient regularization} technique \cite{mishchenko2021regularized,doikov2021gradient},
which allows using the \textit{square} of the Euclidean norm as a regularizer for our model.
Each step of the method becomes easier to compute by employing just one matrix inversion.
Thus, we come to the following scheme:

\begin{algorithm}[h!]
	\caption{Regularized Newton  with Lazy Hessians}
	\label{alg:LazyGradRegNewton}
	\begin{algorithmic}[1]
		\STATE {\bfseries Input:} $\xx_0 \in \R^d$, $m \geq 1$, $L > 0$. Choose $M \geq 0$.
		\FOR{$k=0,1,\dotsc$}
		\STATE Set last snapshot point $\zz_k = \xx_{\pi(k)}$
		\STATE Set regularization parameter $\lambda_k =  \sqrt{M \|\nabla f(\xx_k) \|_{*} }$
		\STATE Compute lazy Newton step: 
		$ \xx_{k + 1}  =  \xx_k - \bigl( \nabla^2 f(\zz_k) + \lambda_k \mat{B} \bigr)^{-1} \nabla f(\xx_k) $
		\ENDFOR
	\end{algorithmic}
\end{algorithm}

The parameter $M$ 
should have the same order and physical interpretation as
the one in Cubic Newton.
We use
\vspace{-1mm}
\beq \label{MGradReg}
\boxed{
	\ba{rcl}
	M & := & 3mL
	\ea
}
\eeq

The global complexity bounds
for Algorithm~\ref{alg:LazyGradRegNewton}
are the same (up to an additive logarithmic term) as those ones
for Algorithm~\ref{alg:LazyCubicNewton}.
However, each iteration of Algorithm~\ref{alg:LazyGradRegNewton}
is much easier to implement since it involves solving just one linear system.

\begin{framedtheorem} \label{ThGradRegGlobal}
Let $M$ be fixed as in \eqref{MGradReg}.
Assume that the gradients for previous iterations $\{ \xx_i \}_{i = 1}^k$
of Algorithm~\ref{alg:LazyGradRegNewton} 
are higher than a desired error level $\varepsilon > 0$:
\beq \label{OldAlgGRBigGrads}
\ba{rcl}
\| \nabla f(\xx_i) \|_{*} & \geq & \varepsilon.
\ea
\eeq
Then, the number of iterations
to reach  accuracy $ \| \nabla f(\xx_{k + 1}) \|_{*} \leq \varepsilon$ is at most
\beq \label{OldAlgGRIters}
\ba{rcl}
\!\!\!\! k & \! \leq \! & 
\cO\Bigl( \frac{\sqrt{mL} (f(\xx_0) - f^\star)}{\varepsilon^{3/2}} 
+ \ln \frac{\| \nabla f(\xx_0) \|_{*}}{\varepsilon} \Bigr).
\ea
\eeq
The total number of Hessian updates $t$ during these iterations is bounded as
\beq \label{OldAlgGRHess}
\ba{rcl}
\!\!\!\! t & \! \leq \! & 
\cO\Bigl( \frac{\sqrt{L} (f(\xx_0) - f^\star)}{\varepsilon^{3/2} \sqrt{m}} 
+ \frac{1}{m} \ln \frac{\| \nabla f(\xx_0) \|_{*}}{\varepsilon} \Bigr).
\ea
\eeq
\end{framedtheorem}

\section{Local Superlinear Convergence}
\label{SectionLocal}

In this section, we discuss local convergence of our second-order schemes
with lazy Hessian updates. We show the fast linear rate with a constant factor during each phase 
and superlinear convergence taking into account the Hessian updates.

Let us assume that our objective is \textit{strongly convex},
thus
$
\nabla^2 f(\xx)  \succeq  \mu \mat{B},  \; \forall \xx \in \R^d,
$
with some parameter $\mu > 0$.
For simplicity, we require that for all points from our space, 
while it is possible to analyze Newton steps assuming strong convexity only in a neighborhood
of the solution. %

\begin{framedtheorem} \label{ThLocal}
Let $M \geq 0$.
Assume that the initial gradient is small enough:
$
\| \nabla f(\xx_0) \|_{*}  \leq  \frac{1}{2}\mathcal{G}_0,
$
where $\mathcal{G}_0 := \frac{\mu^2}{3L + M}$.
Then Algorithm~\ref{alg:LazyCubicNewton} has
a superlinear convergence for the gradient norm, for $k \geq 0$:
$$
\ba{rcl}
\!\!\! \| \nabla f(\xx_k) \|_{*}
& \! \leq \! & 
\mathcal{G}_0 \cdot
\bigl( \frac{1}{2}  \bigr)^{(1 + m)^{\pi(k)}(1 + k \Mymod m) },
\ea
$$
where $\pi(k)$ is defined by \eqref{PiDef}.
\end{framedtheorem}

\BC
Combining both Theorem~\ref{ThGlobCompl} and \ref{ThLocal}, we conclude
that for minimizing a strongly convex function
by Algorithm~\ref{alg:LazyCubicNewton} with regularization parameter 
$M$ given by \eqref{MCubicNewton} and starting from an arbitrary $\xx_0$,
we need 
$$
\ba{rcl}
k & \leq & 
\cO\Bigl( 
\frac{m^2 L^2 (f(\xx_0) - f^\star)}{\mu^3}
\, + \, \frac{1}{\ln(1 + m)} \ln \ln \frac{\mu^2}{mL\varepsilon}
\Bigr)
\ea
$$
lazy steps to achieve $\| \nabla f(\xx_k) \|_{*} \leq \varepsilon$.
\EC

We show also that Algorithm~\ref{alg:LazyGradRegNewton} locally has a
slightly worse but still superlinear convergence rate, which is extremely fast from the practical perspective. For example, for $m=1$, Algorithm~\ref{alg:LazyGradRegNewton} has a superlinear convergence rate of order $3/2$ while Algorithm~\ref{alg:LazyCubicNewton} has order $2$.

\begin{framedtheorem} \label{ThGRLocal}
Let $M \geq 0$.
Assume that the initial gradient is small enough:
$
\| \nabla f(\xx_0) \|_{*}  \leq
\frac{1}{2^2} \mathcal{G}_0,
$
where $\mathcal{G}_0 := \frac{\mu^2}{2^3 (3L + 4M)}$.
Then Algorithm~\ref{alg:LazyGradRegNewton}
has
a superlinear convergence for the gradient norm, for $k \geq 0$:
$$
\ba{rcl}
 \!\!\! \| \nabla f(\xx_k) \|_{*}
& \!\! \leq \!\! &
\mathcal{G}_0 \cdot
\bigl( \frac{1}{2} \bigr)^{2( 1 \, + \, m / 2 )^{\pi(k)} (1 + (k \Mymod m) / 2) },
\ea
$$
where $\pi(k)$ is defined by \eqref{PiDef}.
\end{framedtheorem}

\BC
Combining both Theorem~\ref{ThGradRegGlobal} and \ref{ThGRLocal}, we conclude
that for minimizing a strongly convex function
by Algorithm~\ref{alg:LazyGradRegNewton} with regularization parameter 
$M$ given by \eqref{MGradReg} and starting from an arbitrary $\xx_0$,
we need to do
$$
\ba{rcl}
k & \leq & 
\cO\Bigl( 
\frac{m^2 L^2 (f(\xx_0) - f^\star)}{\mu^3}
\, + \, \ln \frac{mL \|\nabla f(\xx_0) \|_{*}}{\mu^2} \\[10pt]
& &
\, + \, \frac{1}{\ln(1 + m /2)} \ln \ln \frac{\mu^2}{mL\varepsilon}
\Bigr)
\ea
$$
lazy steps to achieve $\| \nabla f(\xx_k) \|_{*} \leq \varepsilon$.
\EC

\section{Practical Implementation}
\label{SectionImplementation}

\subsection{Use of Matrix Factorization}
\label{SubsectionFactorization}

Computing the Hessian once per $m$ iterations, we want to be able to solve efficiently
the corresponding subproblem with this Hessian in each iteration of the method.

For solving the subproblem in Cubic Newton~\eqref{StepDef}, we need to find parameter $r > 0$
that is, in a nondegenerate case, 
the root of the following nonlinear equation (see also Section 5 in \cite{nesterov2006cubic}),
\beq \label{NonlinEq}
\ba{rcl}
\varphi(r) & \Def & \| \ss(r) \| \; - \; r \;\; = \;\; 0,
\ea
\eeq
where
$
\ss(r) \Def 
( \nabla^2 f(\zz) + \frac{M r}{2} \mat{B}  )^{-1} \nabla f(\xx),
$
for $r$ such that the matrix is positive definite.
Applying to \eqref{NonlinEq} the standard bisection method, or the univariate Newton's method
with steps of the form $r^{+} = r - \varphi(r) / \varphi'(r)$, the main operation that we need to perform
is solving a linear equation:
\beq \label{MainLinEquation}
\ba{cl}
& \bigl(\nabla^2 f(\zz) + \tau \mat{B}\bigr) \hh  =  -\nabla f(\xx),
\ea
\eeq
for different values of $\tau > 0$. This is the same type of operation
that we need to do \textit{once} per each step in Algorithm~\ref{alg:LazyGradRegNewton}.

Now, let us assume that we have the following factorization of the Hessian:
\beq \label{Decomp}
\ba{rcl}
\nabla^2 f(\zz) & = & \mat{U} \mat{\Lambda} \mat{U}^\top, 
\quad \text{where} \quad \mat{U} \mat{U}^\top = \mat{B},
\ea
\eeq
and $\mat{\Lambda} \in \R^{d \times d}$ is a \textit{diagonal} or \textit{tridiagonal} matrix.
Thus, $\mat{U}$ is a set of vectors orthogonal with respect to $\mat{B}$.
In case $\mat{B} = \mat{I}$ (identity matrix) and $\mat{\Lambda}$ being diagonal,
\eqref{Decomp} is called \textit{Eigendecomposition}
and implemented in most of the standard Linear Algebra packages.

In general, factorization \eqref{Decomp} can be computed in
$\cO(d^3)$ arithmetic operations. Namely, we can apply a standard 
orthogonalizing decomposition for the following matrix:
$$
\ba{rcl}
\mat{B}^{-1/2} \nabla^2 f(\zz) \mat{B}^{-1/2} & = & 
\mat{V} \mat{\Lambda} \mat{V}^\top, \quad \text{with} \quad
\mat{V}\mat{V}^\top = \mat{I},
\ea
$$
and then set $\mat{U} := \mat{B}^{1/2} \mat{V}$, which gives \eqref{Decomp}.
Note that the solution to \eqref{MainLinEquation} can expressed as
$$
\ba{rcl}
\hh & = & - \mat{U}^{-\top} \bigl( \mat{\Lambda} + \tau \mat{I} \bigr)^{-1} \mat{U}^{-1} \nabla f(\xx),
\ea
$$
and it is computable in $\cO(d^2)$ arithmetic operations for any given 
$\tau > 0$ and $\nabla f(\xx)$.
The use of tridiagonal decomposition can be more efficient in practice.
Indeed, inversion of the matrix $\mat{\Lambda} + \tau \mat{I}$
would still cost $\cO(d)$ operations for tridiagonal $\mat{\Lambda}$, while 
it requires less computational resources and less floating point precision 
to compute such decomposition.
In practice, it is also important to leverage a structure
of the Hessian (e.g. when matrices are \textit{sparse}), which can further improve the arithmetic cost of each step.

\subsection{Adaptive Search}
\label{SubsectionAdaptive}

To obtain our results we needed to pick the regularization parameter $M $ proportional to $mL$,
where $m \geq 1$ is the length of one phase
and $L > 0$ is the Lipschitz constant of the Hessian. One drawback of this choice is that we actually need to know the Lipschitz constant,
which is not always the case in practice. Moreover, with a constant regularization,
the methods become conservative, preventing the possibility of big steps. 
At the same time, from the local perspective, the best quadratic approximation of the objective
is the pure second-order Taylor polynomial. So, being in a neighborhood of the solution,
the best strategy is to tend $M$ to zero which gives the fastest local superlinear rates
(see Section~\ref{SectionLocal}).

The use of adaptive search in second-order schemes has been studied for several years
\cite{nesterov2006cubic,cartis2011adaptive1,cartis2011adaptive2,grapiglia2017regularized,grapiglia2019accelerated,doikov2021minimizing}.
It is well known that such schemes have a very good practical performance,
while the adaptive search makes the methods to be also \textit{universal} \cite{grapiglia2017regularized,doikov2021minimizing},
that is to adapt to the best H\"older degree of the Hessian,
or even \textit{super-universal} \cite{doikov2022super} which 
chooses automatically the H\"older degree of either the second or third derivative of the objective, 
working properly
on a wide range of problem classes with the best global complexity guarantees.

We propose Algorithm~\ref{alg:AdaptiveCubicNewton} which changes the value of $M$ adaptively and checks the functional progress after $m$ steps to validate its choice. Importantly, the knowledge of the Lipschitz constant $L$ is not needed.

\begin{algorithm}[h!]
	\caption{Adaptive Cubic Newton with Lazy Hessians}
	\label{alg:AdaptiveCubicNewton}
	\begin{algorithmic}[1]
		\STATE {\bfseries Input:} $\xx_0 \in \R^d$, $m \geq 1$. Fix some $M_0 > 0$.
		\FOR{$t=0,1,\dotsc$}
		\STATE Compute snapshot Hessian $\nabla^2 f( \xx_{tm} )$
		\REPEAT
		\STATE Update $M_t = 2 \cdot M_t$
		\FOR{$i=1,\ldots,m$}
		\STATE Lazy cubic step $\xx_{tm + i} = T_{M_t}(\xx_{tm + i - 1}, \xx_{tm})$
		\ENDFOR
		\UNTIL {\small $f(\xx_{tm})  -  f(\xx_{tm + m})  \geq 
			\frac{1}{\sqrt{M_t}} \sum\limits_{i = 1}^m \! 
			\| \nabla f(\xx_{tm + i})   \|_{*}^{3/2}$}
		\STATE Set $M_{t + 1} = \frac{1}{4} \cdot M_t$
		\ENDFOR
	\end{algorithmic}
\end{algorithm}

 $M_0$ is an \textit{initial guess} for the regularization constant,
which can be further both increased or decreased dynamically. 

\begin{framedtheorem} \label{ThAdaptiveCubicGloblal}
Let the gradients for previous iterates $\{ \xx_i \}_{i = 1}^{tm}$ 
during $t$ phases of Algorithm~\ref{alg:AdaptiveCubicNewton}
are higher than a desired error level $\varepsilon > 0$:
\beq \label{OldAdAlgGLB}
\ba{rcl}
\| \nabla f(\xx_i) \|_{*} & \geq & \varepsilon.
\ea
\eeq
Then, the number of phases
to reach accuracy $\| \nabla f(\xx_{tm + 1}) \|_{*} \leq \varepsilon$
is at most
\beq \label{OldAdAlgPhases}
\ba{rcl}
t & \leq & \cO\Bigl( 
\sqrt{
	\max\bigl\{ \frac{M_0}{m}, L \bigr\} } \cdot
\frac{f(\xx_0) - f^\star}{\varepsilon^{3/2}\sqrt{m}}  \Bigr).
\ea
\eeq
The total number $N$ of gradient calls during these phases is bounded as
\beq \label{AdAlgTotal}
\ba{rcl}
N & \leq & 
2tm
\, + \,
\max\{1, \log_2 \frac{2^9 3^5 mL}{M_0}  \} m.
\ea
\eeq
\end{framedtheorem}

\BR
According to Theorem \ref{ThAdaptiveCubicGloblal}, 
the global complexity of the adaptive algorithm
is the same as the one given by Theorem~\ref{ThGlobCompl} 
for the method with a fixed regularization constant.
Due to \eqref{AdAlgTotal}, the average number
of tries of different $M_t$ per one phase is only \underline{two}.
\ER

Clearly, It is also possible to incorporate a similar adaptive search into Algorithm~\ref{alg:LazyGradRegNewton}
for the convex case (see Appendix~\ref{ASCV}).

\section{Experiments}
\label{SectionExperiments}

We demonstrate an illustrative numerical experiment on the performance of the proposed second-order methods
with lazy Hessian updates.
We consider the following convex minimization problem
with the \textit{Soft Maximum} objective (log-sum-exp):
\beq \label{LogSumExp}
\ba{rcl}
\min\limits_{\xx \in \R^d} f(\xx) & := & 
\mu \ln\biggl(\,  
\sum\limits_{i = 1}^n \exp\Bigl(  \frac{\la \aa_i, \xx \ra - b_i}{\mu} \Bigr) 
\,\biggr) \\
\\
& \approx & 
\max\limits_{1 \leq i \leq n} \Bigl[ \la \aa_i, \xx \ra - b_i \Bigr].
\ea
\eeq
The problems of this type are important in applications with minimax 
strategies for matrix games and for training $\ell_{\infty}$-regression \cite{nesterov2005smooth,bullins2020highly}.

To generate the data, we sample randomly
the vectors $\bar{\aa}_1, \ldots, \bar{\aa}_n \in \R^d$ and $\bb \in \R^n$
with elements from the uniform distribution on $[-1 , 1]$.
Then, we build an auxiliary objective $\bar{f}$
of the form \eqref{LogSumExp} with these vectors and set
$\aa_{i} := \bar{\aa}_{i} - \nabla \bar{f}(\0)$. This ensures 
the optimum is at the origin since $\nabla f(\0) = \0$.
The starting point is $\xx_0 = (1, \ldots, 1)$.

For the primal norm~\eqref{EuclideanNorm}, we use the matrix
\beq \label{BExample}
\ba{rcl}
\mat{B} & := & \sum\limits_{i = 1}^n \aa_i \aa_i^\top  +  \delta \mat{I} \;\; \succ \;\; 0,
\ea
\eeq
where $\delta > 0$ is a small perturbation parameter to ensure positive definiteness.
Then, the corresponding Lipschitz constant of the Hessian is bounded by
(see, e.g. \citep[Example 1.3.5]{doikov2021new}):
$
L   = 2 / \mu^2,
$
where $\mu > 0$ is a smoothing parameter.

\begin{figure}[h!]
	\centering
	\hspace*{-5pt}
	\includegraphics[width=0.49\textwidth ]{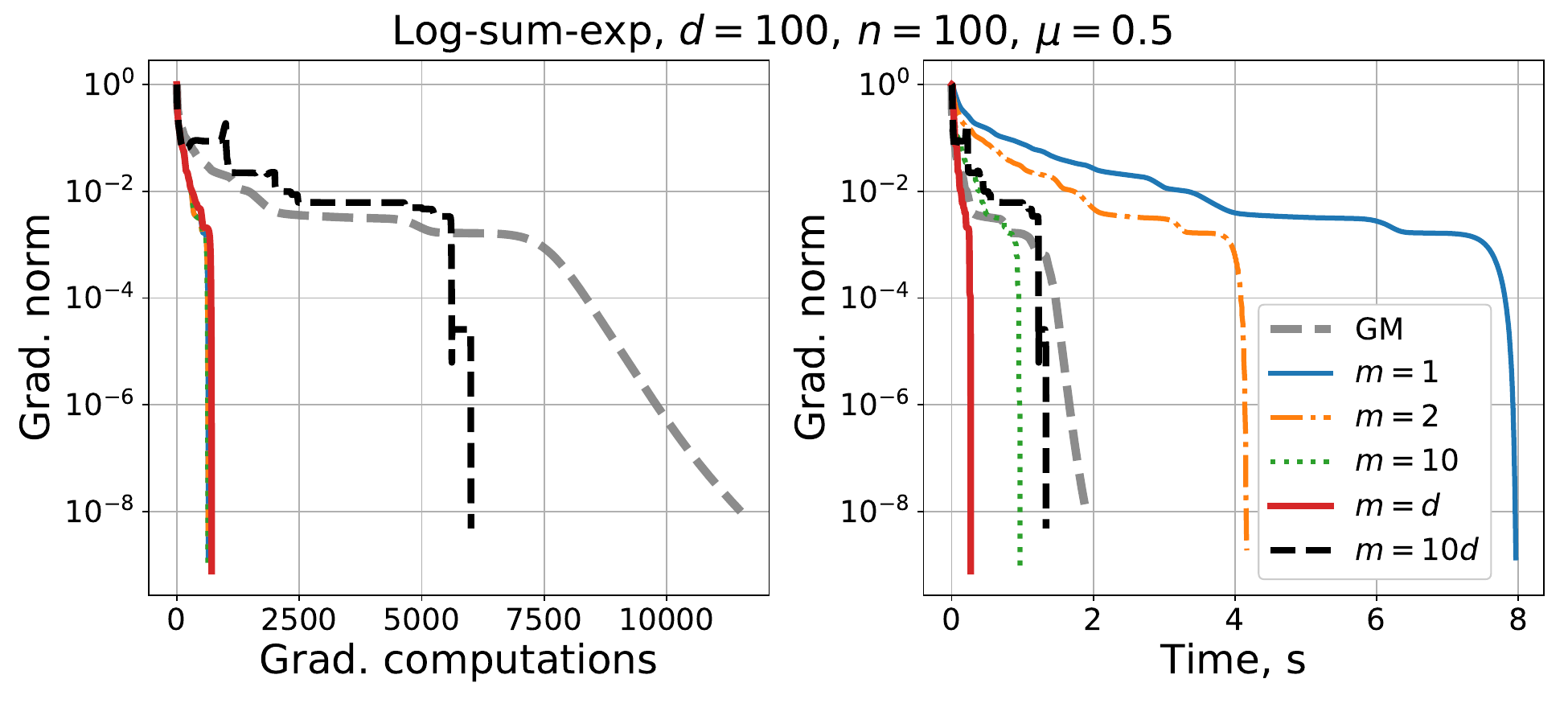} \\
	\hspace*{-5pt}
	\includegraphics[width=0.49\textwidth ]{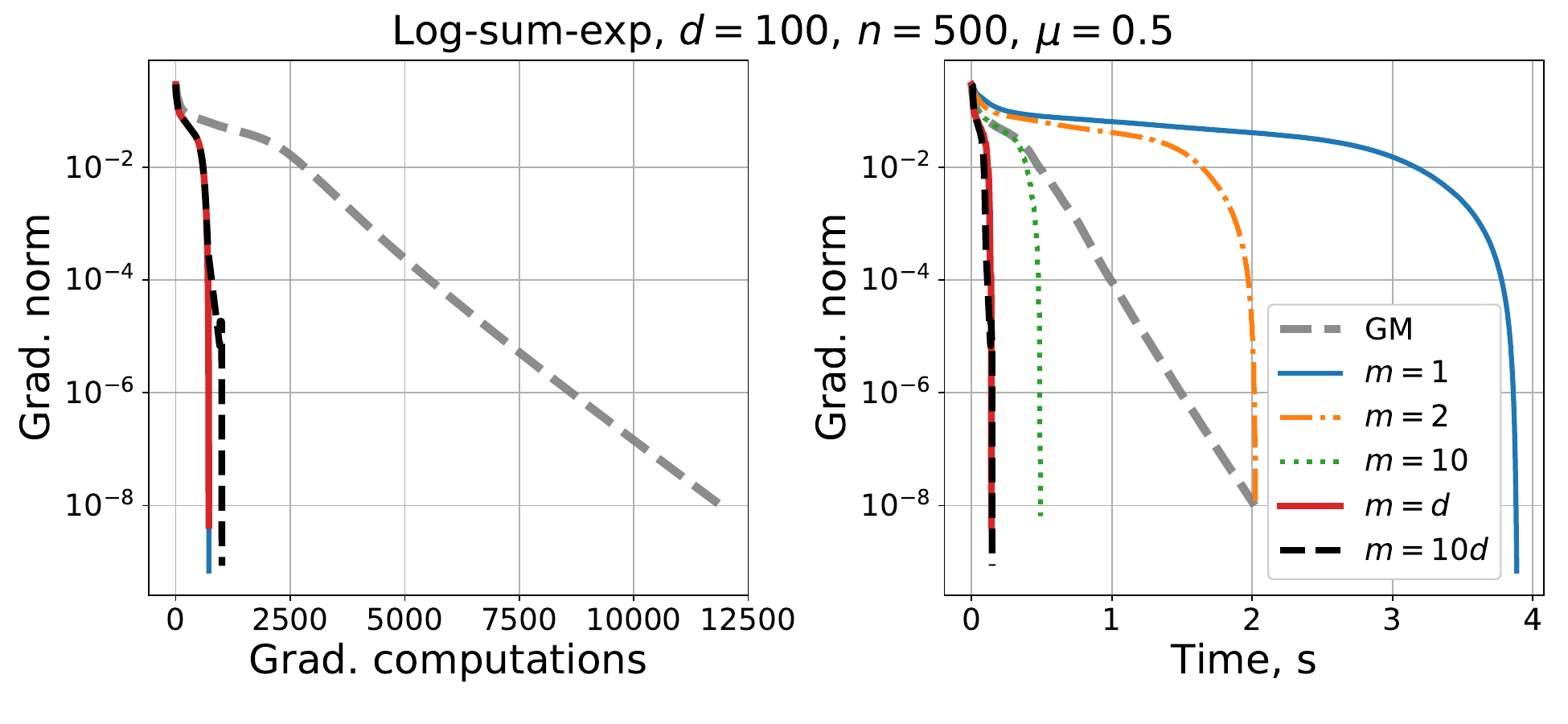}
	\caption{\small Gradient norm depending on the number of computed gradients and on running time, when varying the frequency $m$ of updating the Hessian in Algorithm~\ref{alg:LazyGradRegNewton}.}
	\label{fig:FixedM}
\end{figure}

Since the problem is convex, we can apply 
Newton's method with gradient regularization (Algorithm~\ref{alg:LazyGradRegNewton}).
In Figure~\ref{fig:FixedM}, we compare different values of parameter $m$ that is the frequency of updating the Hessian.
The regularization parameter is fixed as $M := 1$. We also show the performance of the Gradient Method as a standard baseline.

\smallskip

We see that increasing the parameter $m$, thus reusing old Hessians for more of the inner steps, we significantly improve the overall performance of the method
in terms of the total computational time. The best frequency is $m := d$ which confirms our theory.

\section{Discussion}
\label{SectionDiscussion}

\paragraph{Conclusions.}
In this work, we have developed new second-order algorithms
with \textit{lazy Hessian updates}
for solving general possibly non-convex optimization problems. 
We show that it can be very efficient to reuse the previously computed Hessian 
for several iterations of the method,
instead of updating the Hessian each step.

Indeed, in general, the cost of computing the Hessian matrix is $d$ times more expensive than
computing one gradient.
At the same time, it is intuitively clear that even inexact second-order information 
from the past should significantly help the method  in dealing with ill-conditioning of the problem.
In our analysis, we show that this intuition truly works.

By using cubic regularization and gradient regularization techniques, 
we establish fast global and local rates for our second-order methods with lazy Hessian updates.
We show that the optimal strategy for updating the Hessian is once per~$d$ iterations,
which gives a provable improvement of the total arithmetic complexity by a factor of $\sqrt{d}$.
Our approach also works with classical Newton steps
(without regularization), achieving a local quadratic convergence.

Note that it is possible to 
	extend our results onto the \textit{composite} convex optimization problems, i.e. 
	minimizing the sum $f(\xx) + \psi(\xx)$, where $f(\cdot)$ is convex and smooth,
	while $\psi(\cdot)$ is
	a \textit{simple} closed convex proper function but
	not necessarily differentiable (see Appendix~\ref{SectionCCase}).

\paragraph{Directions for Future Work.} One important direction for further research can be
a study of the problems with a \textit{specific Hessian structure}
(e.g. sparsity or a certain spectral clustering). 
Then, we may need to have different schedules for updating the Hessian matrix,
while it can also help in solving each step more efficiently.

Another interesting question is to make connections between our approach and
classical \textit{Quasi-Newton} methods \cite{nocedal2006numerical}, 
which gradually update the approximation of the Hessian after each step.
Recently discovered non-asymptotic complexity bounds \cite{rodomanov2021new, rodomanov2022quasi} 
for Quasi-Newton methods
may be especially useful for reaching this goal.

We also think that it is possible to generalize our analysis to high-order optimization schemes \cite{nesterov2021implementable} as well.
Note that the main advantage of our schemes is that we can reuse the precomputed factorization of the Hessian and thus we can perform several steps with the same Hessian in an efficient manner. At the same time, it is not clear up to now, how to compute the tensor step efficiently by reusing an old high-order tensor. This would require using some advanced tensor decomposition techniques.

Finally, it seems to be very important to study the effect of lazy Hessian updates for convex optimization
in more details. In our analysis, we studied only a general non-convex convergence in terms of the gradient norm.
Another common accuracy measure in the convex case is the \textit{functional residual}.
Thus, it could be possible to prove some better convergence rates using this measure,
as well as considering \textit{accelerated}  \cite{nesterov2018lectures}
and \textit{super-universal} \cite{doikov2022super} second-order schemes.

We keep these questions for further investigation.

\section*{Acknowledgements}

This work was supported by the Swiss State Secretariat for Education, Research and Innovation (SERI)
under contract number 22.00133.

\bibliography{bibliography}	 
\bibliographystyle{icml2023}

\clearpage
\onecolumn
\appendix

\icmltitle{Supplementary Material}

We consider the problem~\eqref{MainProblem} : \beq
\ba{c}
\min\limits_{\xx \in \R^d} f(\xx),
\ea
\eeq

We assume that the Hessian of $f$ is Lipschitz continuous, with some constant $L > 0$:
\beq \label{LipHess}
\ba{rcl}
\| \nabla^2 f(\xx) - \nabla^2 f(\yy) \| & \leq & L \|\xx - \yy\|, \qquad \forall \xx, \yy \in \R^d.
\ea
\eeq
The immediate consequences are the following inequalities, valid for all $\xx, \yy \in \R^d$:
\beq \label{LipHessGrad}
\ba{rcl}
\| \nabla f(\yy) - \nabla f(\xx) - \nabla^2 f(\xx)(\yy - \xx)\|_{*} & \leq & \frac{L}{2}\|\xx - \yy\|^2,
\ea
\eeq
\beq \label{LipHessFunc}
\ba{rcl}
| f(\yy) - f(\xx) - \la \nabla f(\xx), \yy - \xx \ra - \frac{1}{2} \la \nabla^2 f(\xx)(\yy - \xx), \yy - \xx \ra |
& \leq & \frac{L}{6}\|\yy - \xx\|^3.
\ea
\eeq

In our analysis, we frequently use  the standard Young's inequality for product:
$$
\ba{rcl}
ab & \leq & \frac{a^p}{p} + \frac{b^{q}}{q}, \qquad \forall a, b \geq 0,
\ea
$$
which is valid for any $p, q > 1$ such that $\frac{1}{p} + \frac{1}{q} = 1$. 
\section{Proofs for Section~\ref{SectionLazyCubic}}

Our lazy Newton steps,
 allow the gradient and Hessian to be computed at different points $\xx, \zz \in \R^d$.
We use cubic regularization of our model with a parameter $M > 0$:
\beq \label{App:StepDef}
\ba{rcl}
\mat{T}_{M}(\xx, \zz) & \in & \Argmin\limits_{\yy \in \R^d}
\Bigl\{ 
\la \nabla f(\xx), \yy - \xx \ra + \frac{1}{2} \la \nabla^2 f(\zz)(\yy - \xx), \yy - \xx \ra
+ \frac{M}{6}\|\yy - \xx\|^3
\Bigr\}.
\ea
\eeq
For $\zz := \xx$ this is the 
iteration of Cubic Newton~\cite{nesterov2006cubic}. 
Thus, in our scheme, we can reuse the Hessian from a previous step $\zz$ without recomputing it,
which significantly reduces the overall iteration cost.  

Our definition implies that the point $\mat{T} = \mat{T}_{M}(\xx, \zz)$
is a \textit{global minimum} of the cubically regularized model,
which is generally non-convex. However, it turns out that we can compute this point efficiently
by using standard techniques developed initially for trust-region methods~\cite{conn2000trust}. 

Let us denote $r \Def \|\mat{T} - \xx\|$. 
The solution to the subproblem~\eqref{App:StepDef} satisfies the following stationarity conditions (see \cite{nesterov2006cubic}):
\beq \label{CubicStat}
\ba{rcl}
\nabla f(\xx) + \nabla^2 f(\zz)(\mat{T} - \xx) + \frac{Mr}{2} \mat{B}(\mat{T} - \xx) & = & \0,
\ea
\eeq
and
\beq \label{CubicPSD}
\ba{rcl}
\nabla^2 f(\zz) + \frac{Mr}{2} \mat{B} & \succeq & 0.
\ea
\eeq

Thus, in the non-degenerate case, one step can be represented in the following form: 
\beq \label{CubicStep}
\ba{rcl}
\mat{T} & = & \xx - \bigl( \nabla^2 f(\zz) + \frac{Mr}{2} \mat{B} \bigr)^{-1} \nabla f(\xx),
\ea
\eeq
and the value $r > 0$ can be found by solving the corresponding univariate nonlinear equation
(see \cite{nesterov2006cubic}[Section 5]).
It can be done very efficiently from a precomputed \textit{eigenvalue} or the \textit{tridiagonal}
decomposition of the Hessian. These decompositions typically require $\cO(d^3)$ arithmetic operations,
which is of similar cost as for matrix inversion in the classical Newton step.
We discuss the computation of the new iterate $\mat{T}_{M}(\xx, \zz)$ in more detail in Section~\ref{SubsectionFactorization}.

Now, let us express the progress achieved by the proposed step \eqref{App:StepDef}. 
Firstly, we can relate the length of the step and the gradient at the new point as follows:
$$
\ba{rcl}
\frac{Lr^2}{2} & \overset{\eqref{LipHessGrad}}{\geq} &
\| \nabla f(\mat{T}) - \nabla f(\xx) - \nabla^2 f(\xx)(\mat{T} - \xx) \|_{*} \\
\\
& \overset{\eqref{CubicStat}}{=} &
\| \nabla f(\mat{T}) + \frac{Mr }{2} \mat{B}(\mat{T} - \xx) + (\nabla^2 f(\zz) - \nabla^2 f(\xx))(\mat{T} - \xx) \|_{*} \\
\\
& \geq & 
\| \nabla f(\mat{T}) \|_{*} - \frac{Mr^2}{2} - r \| \nabla^2 f(\zz) - \nabla^2 f(\xx) \|.
\ea
$$
Rearranging terms and using Lipschitz continuity of the Hessian, we get
$$
\ba{rcl}
\| \nabla f(\mat{T}) \|_{*} & \leq & 
\frac{M + L}{2} r^2
+ r \| \nabla^2 f(\zz) - \nabla^2 f(\xx) \| \\
\\
& \overset{\eqref{LipHess}}{\leq} &
\frac{M + L}{2} r^2 + L r \|\zz - \xx\|.
\ea
$$
Thus, we have established
\BL \label{LemmaCubicGrad}
For any $M > 0$, it holds that
\beq \label{NewGradBound}
\ba{rcl}
\| \nabla f(\mat{T}) \|_{*} & \leq & \frac{M + L}{2} r^2 + L r\|\zz - \xx\|.
\ea
\eeq
\EL

Secondly, we have the following progress in terms of the objective function.

\BL \label{LemmaCubicFunc}
For any $M > 0$, it holds that
\beq \label{NewFuncProgress}
\ba{rcl}
f(\xx) - f(\mat{T}) & \geq & \frac{5M - 4L}{24} r^3
- \frac{32L^3}{3M^2}\|\zz - \xx\|^3.
\ea
\eeq
\EL
\proof
Indeed, multiplying \eqref{CubicStat} by the vector $\mat{T} - \xx$, we get
\beq \label{GradProd}
\ba{rcl}
\la \nabla f(\xx), \mat{T} - \xx \ra & = &
- \la \nabla^2 f(\zz)(\mat{T} - \xx), \mat{T} - \xx \ra - \frac{M}{2} r^3 \\
\\
& \overset{\eqref{CubicPSD}}{\leq} &
- \frac{1}{2} \la \nabla^2 f(\zz)(\mat{T} - \xx), \mat{T} - \xx \ra - \frac{M}{4}r^3.
\ea
\eeq
Hence, using the global upper bound on the objective, we obtain
$$
\ba{rcl}
f(\mat{T}) & \overset{\eqref{LipHessFunc}}{\leq} &
f(\xx) + \la \nabla f(\xx), \mat{T} - \xx \ra + \frac{1}{2}\la \nabla^2 f(\xx)(\mat{T} - \xx), \mat{T} - \xx \ra 
+ \frac{L}{6} r^3 \\
\\
& \overset{\eqref{GradProd}}{\leq} &
f(\xx) - \frac{M}{4}r^3 + \frac{L}{6}r^3 + \frac{1}{2} \la ( \nabla^2 f(\xx) - \nabla^2 f(\zz) ) (\mat{T} - \xx), \mat{T} - \xx \ra \\
\\
& \overset{\eqref{LipHess}}{\leq} &
f(\xx) - \frac{M}{4}r^3 + \frac{L}{6}r^3 + \frac{L}{2} r^2 \|\zz - \xx\|.
\ea
$$
Applying Young's inequality for the last term,
$$
\ba{rcl}
\frac{L}{2} r^2 \|\zz - \xx\| 
& = &
\bigl( \frac{M^{2/3}}{2 \cdot 32^{1/3}} r^2 \bigr) 
\cdot \bigl(  \frac{32^{1/3} L}{M^{2/3}}  \|\zz - \xx\| \bigr)
\;\; \leq \;\; \frac{M}{24} r^3 + \frac{32 L^3}{3 M^2} \|\zz - \xx\|^3,
\ea
$$
we obtain \eqref{NewFuncProgress}. \qed

We can also bound the smallest eigenvalue of the new Hessian, which will be crucial to understand the behaviour for non-convex objectives. Indeed, 
using Lipschitz continuity of the Hessian and the triangle inequality, we have:
\beq \label{HessEig}
\ba{rcl}
\nabla^2 f(\mat{T}) & \succeq &
\nabla^2 f(\zz) - L\|\zz - \mat{T}\| \mat{B} \\
\\
& \succeq &
\nabla^2 f(\zz) - (L\|\zz - \xx\|  + Lr) \mat{B} \\
\\
& \overset{\eqref{CubicPSD}}{\succeq} &
- (\frac{M}{2}r + L\|\zz - \xx\|  + Lr)\mat{B}.
\ea
\eeq

Let us denote the following quantity, for any $y \in \R^d$:
\beq \label{XiDef}
\ba{rcl}
\xi(\yy) & \Def & \Bigl[  -\lambda_{\min}\bigl( 
\mat{B}^{-1/2} \nabla^2 f(\yy) \mat{B}^{-1/2}
\bigr) \Bigr]_{+},
\ea
\eeq
where $[t]_{+} \Def \max\{ t, 0 \}$ denotes the positive part,
and $\lambda_{\min}(\cdot)$ is the smallest eigenvalue of a symmetric matrix.

If $\nabla^2 f(\yy) \succeq 0$ for a certain $\yy \in \R^d$, then $\xi(\yy) = 0$.
Otherwise, $\xi(\yy)$ shows how big (in absolute value) 
the smallest eigenvalue of the Hessian is with respect to a fixed matrix $\mat{B} \succ 0$.

Thus, due to \eqref{HessEig}, we can control the quantity $\xi$ as follows:

\BL \label{LemmaCubicHess}
For any $M > 0$, it holds
\beq \label{NewHessBound}
\ba{rcl}
\xi(\mat{T}) & \leq &  \frac{M + 2L}{2} r + L \|\zz - \xx\|.
\ea
\eeq
\EL

Now, we can combine \eqref{NewGradBound}
and \eqref{NewHessBound} with \eqref{NewFuncProgress}
to express the functional progress of one step using the new gradient norm
and the smallest eigenvalue of the new Hessian.

\BT \label{ThProgressApp}
Let $M \geq L$. Then, for one cubic step~\eqref{StepDef}, we have the one step progress bound \eqref{OneStepProgress}, i.e.,
\beq 
\ba{rcl}
f(\xx) - f(\mat{T}) & \geq & 
\max\Bigl\{ \frac{1}{648M^2}\xi(\mat{T})^3, \;
\frac{1}{72 \sqrt{2M}} \| \nabla f(\mat{T}) \|_{*}^{3 / 2}
\Bigr\}
+ \frac{M}{48} r^3
- \frac{11L^3}{M^2}\|\zz - \xx\|^3.
\ea
\eeq
\ET
\proof
Indeed, using convexity of the function $t \mapsto t^3$ for $t \geq 0$, we obtain
$$
\ba{rcl}
\xi(\mat{T})^3 & \overset{\eqref{NewHessBound}}{\leq} &
\Bigl(  \frac{3}{2} Mr  + L \|\zz - \xx\| \Bigr)^3
\;\; \leq \;\;
\frac{27}{2} M^3 r^3 + 4 L^3 \|\zz - \xx\|^3.
\ea
$$
Hence,
$$
\ba{rcl}
\frac{Mr^3}{2} & \geq & 
\frac{1}{27 M^2} \xi(\mat{T})^3
- \frac{4L^3}{27M^2} \|\zz - \xx\|^3,
\ea
$$
which gives the first part of the maximum when substituting into \eqref{NewFuncProgress}.

Then, using convexity of the function $t \mapsto t^{3/2}$ for $t \geq 0$ and Young's inequality,
we get
$$
\ba{rcl}
\| \nabla f(\mat{T}) \|_{*}^{3/2} & \overset{\eqref{NewGradBound}}{\leq} &
\Bigl(  M r^2 + Lr \|\zz - \xx\| \Bigr)^{3/2} 
\;\; \leq \;\; \sqrt{2} M^{3/2} r^3
+ \sqrt{2} (Lr \|\zz - \xx\|)^{3/2} \\
\\
& \leq & 
\sqrt{2} M^{3/2} r^3 + \frac{\sqrt{2} M^{3/2}}{2} r^3 
+ \frac{\sqrt{2} L^{3}}{2 M^{3/2}} \|\zz - \xx\|^3
\;\; = \;\;
\frac{3 \sqrt{2} M^{3/2}}{2} r^3 + \frac{L^3}{\sqrt{2} M^{3/2}} \|\zz - \xx\|^3.
\ea
$$
Therefore,
$$
\ba{rcl}
\frac{M r^3}{2} & \geq & 
\frac{1}{3 \sqrt{2 M}} \| \nabla f(\mat{T}) \|_{*}^{3/2}
- \frac{L^3}{6 M^2} \|\zz - \xx\|^3.
\ea
$$
Substituting this bound into~\eqref{NewFuncProgress} completes the proof.
\qed

\section{Proofs for Section~\ref{SectionGlobalCubic}}
\label{AppendixSectionGlobal}

Let us denote the corresponding step lengths by $r_{k + 1} := \|\xx_{k + 1} - \xx_{k}\|$.
Then, according to Theorem~\ref{ThProgress}, for each $0 \leq k \leq m - 1$, we have the following
progress in terms of the objective function:
$$
\ba{rcl}
f(\xx_k) - f(\xx_{k + 1}) & \overset{\eqref{OneStepProgress}}{\geq} &
\max\Bigl\{
\frac{1}{648M^2}\xi(\xx_{k + 1})^3, \,
\frac{1}{72 \sqrt{2M}} \| \nabla f(\xx_{k + 1}) \|_{*}^{3/2}
\Bigr\}
+ \frac{M}{48} r_{k + 1}^3 - \frac{11L^3}{M^2}\|\xx_0 - \xx_{k}\|^3.
\ea
$$
Telescoping this bound for different $k$, and using triangle inequality
for the last negative term,
$$
\ba{rcl}
\|\xx_0 - \xx_k\| 
& \leq & 
\sum\limits_{i = 1}^{k} r_i, 
\ea
$$
we get
\beq \label{OnePhaseBound}
\ba{cl}
& f(\xx_0) - f(\xx_{m}) \\
\\
& \geq \;\; 
\sum\limits_{k = 1}^m 
\max\Bigl\{
\frac{1}{648M^2}\xi(\xx_k)^3, \,
\frac{1}{72\sqrt{2M}} \sum\limits_{k = 1}^m \| \nabla f(\xx_k) \|_{*}^{3/2} \Bigr\}
+ \frac{M}{48} \sum\limits_{k = 1}^m r_k^3
- \frac{11L^3}{M^2} \sum\limits_{k = 1}^{m - 1} \Bigl(  \sum\limits_{i = 1}^{k} r_i  \Bigr)^3.
\ea
\eeq

It remains to use the following simple technical lemma.
\BL
For any sequence of positive numbers $\{ r_k \}_{k \geq 1}$, it holds for any $m \geq 1$:
\beq \label{TechBound}
\ba{rcl}
\sum\limits_{k = 1}^{m - 1} \Bigl(  \sum\limits_{i = 1}^{k} r_i  \Bigr)^3 
& \leq & \frac{m^3}{3} \sum\limits_{k = 1}^{m - 1} r_k^3
\ea
\eeq
\EL
\proof
We prove \eqref{TechBound} by induction. It is obviously true for $m = 1$, which is our base.

Assume that it holds for some arbitrary $m \geq 1$. Then
\beq \label{Induct1}
\ba{rcl}
\sum\limits_{k = 1}^{m} \Bigl(  \sum\limits_{i = 1}^{k} r_i  \Bigr)^3 
& = &
\sum\limits_{k = 1}^{m - 1} \Bigl(  \sum\limits_{i = 1}^{k} r_i  \Bigr)^3 
+ \Bigl(\sum\limits_{i = 1}^m r_i \Bigr)^3
\;\; \overset{\eqref{TechBound}}{\leq} \;\;
\frac{m^3}{3} \sum\limits_{k = 1}^{m - 1} r_k^3
+ \Bigl(\sum\limits_{i = 1}^m r_i \Bigr)^3.
\ea
\eeq
Applying Jensen's inequality for
convex function $t \mapsto t^3$, $t \geq 0$, we have
a bound for the second term:
\beq \label{Induct2}
\ba{rcl}
\Bigl( \sum\limits_{i = 1}^m r_i \Bigr)^3
\;\; = \;\;
m^3 \Bigl( \frac{1}{m} \sum\limits_{i = 1}^m r_i \Bigr)^3
& \leq & m^2 \sum\limits_{i = 1}^m r_i^3.
\ea
\eeq
Therefore,
$$
\ba{rcl}
\sum\limits_{k = 1}^{m} \Bigl(  \sum\limits_{i = 1}^{k} r_i  \Bigr)^3 
& \overset{\eqref{Induct1}, \eqref{Induct2}}{\leq} &
\Bigl(  \frac{m^3}{3} + m^2 \Bigr) \sum\limits_{k = 1}^m r_k^3
\;\; \leq \;\;
\frac{(m + 1)^3}{3} \sum\limits_{k = 1}^m r_k^3. \QF
\ea
$$

Using this bound, we ensure the right-hand side of \eqref{OnePhaseBound}
to be positive. For that, we need to use an increased value of the regularization parameter.
Thus, we prove the following guarantee.
\BC
Let $M \geq 6 m L$. Then,
\beq \label{Telescoped}
\ba{rcl}
f(\xx_0) - f(\xx_m) & \geq &
\sum\limits_{k = 1}^m
\max\Bigl\{
\frac{1}{648M^2} \xi(\xx_k)^3, \,
\frac{1}{72 \sqrt{2M}}  \| \nabla f(\xx_k)\|_{*}^{3/2}\Bigr\}.
\ea
\eeq
\EC

\BT \label{ThGlobComplApp}
Let the objective function be bounded from below: $\inf\limits_{\xx} f(\xx) \geq f^\star$.
Let the regularization parameter $M$ be fixed as in \eqref{MCubicNewton}. Assume that the gradients
for previous iterations $\{ \xx_i \}_{i = 1}^k$ are higher than a desired error level $\varepsilon>0$:\vspace{-1mm}
\beq \label{GradBounded}
\ba{rcl}
\| \nabla f(\xx_i) \|_{*} & \geq & \varepsilon.
\ea
\eeq
Then, the number of iterations of Algorithm~\ref{alg:LazyCubicNewton} to reach accuracy $\| \nabla f(\xx_{k + 1}) \|_{*} \leq \varepsilon$ is at most
\beq \label{BoundAllItersAPP}
\ba{rcl}
k & \leq & \cO\Bigl( \frac{\sqrt{m L} (f(\xx_0) - f^\star) }{\varepsilon^{3/2}} \Bigr),
\ea
\eeq
where $\cO(\cdot)$ hides some absolute numerical constants.
The total number of Hessian updates $t$ during these iterations is bounded as \eqref{BoundAllHess}, i.e.,
\beq
\ba{rcl}
t & \leq & \cO\Bigl( \frac{\sqrt{L} (f(\xx_0) - f^\star)}{\varepsilon^{3/2} \sqrt{m}} \Bigr).
\ea
\eeq
For the minimal eigenvalues of all Hessians, it holds that
\beq \label{App:HessEigenBound}
\ba{rcl}
\min\limits_{1 \leq i \leq k}
\Bigl[-\lambda_{\min}\bigl( \mat{B}^{-1/2} \nabla^2 f(\xx_i) \mat{B}^{-1/2} \bigr)\Bigr]_{+}
& \leq &
\Bigl(\frac{648M^2(f(\xx_0) - f^\star)}{k} \Bigr)^{1/3}.
\ea
\eeq
\ET
\proof
Without loss of generality, we assume that $k$ is a multiple of $m$: $k = t m$.
Then, for the $i$th ($1 \leq i \leq t$) phase of the method, we have
\beq \label{OnePhaseBoundTh}
\ba{rcl}
f(\xx_{m(i - 1)}) - f(\xx_{mi}) & \overset{\eqref{Telescoped}, \eqref{GradBounded}}{\geq} &
\frac{m}{72\sqrt{2M}} \varepsilon^{3/2}
\;\; = \;\;
\frac{\sqrt{m}}{144\sqrt{3L}} \varepsilon^{3/2}.
\ea
\eeq
Telescoping this bound for all phases, we obtain\vspace{-1mm}
$$
\ba{rcl}
f(\xx_0) - f^\star & \geq & f(\xx_0) - f(\xx_k)
\;\; \overset{\eqref{OnePhaseBoundTh}}{\geq} \;\;
\frac{t \sqrt{m}}{144\sqrt{3L}} \varepsilon^{3/2}.
\ea
$$
This gives the claimed bound \eqref{BoundAllHess} for $t$.
Taking into account that $k = tm$,
\eqref{BoundAllItersAPP} follows immediately.
Telescoping the bound~\eqref{Telescoped} for the smallest Hessian eigenvalue $\xi(\cdot)$, we obtain
\beq \label{EigenTelescoped}
\ba{rcl}
f(\xx_0) - f^\star & \geq & f(\xx_0) - f(\xx_k) \;\; \geq \;\; 
\frac{k}{648M^2} \min\limits_{1 \leq i \leq k} \xi(\xx_i)^3.
\ea
\eeq
Hence, \vspace{-1mm}
$$
\ba{rcl}
\min\limits_{1 \leq i \leq k}
\Bigl[-\lambda_{\min}\bigl( \mat{B}^{-1/2} \nabla^2 f(\xx_i) \mat{B}^{-1/2} \bigr)\Bigr]_{+}
& \equiv &
\min\limits_{1 \leq i \leq k} \xi(\xx_i) 
\;\; \overset{\eqref{EigenTelescoped}}{\leq} \;\;
\Bigl(\frac{648 M^2(f(\xx_0) - f^\star)}{k} \Bigr)^{1/3},
\ea
$$
which completes the proof.
\qed

\section{Proofs for Section~\ref{SectionConvex}}
Here we assume all the Hessians are positive semidefinite:
\beq \label{Convexity}
\ba{rcl}
\nabla^2 f(\xx) & \succeq & 0, \qquad \forall \xx\in \R^d.
\ea
\eeq

For a fixed $\zz \in \R^d$, let us denote one 
regularized lazy Newton step from point $\xx \in \R^d$ by
\beq \label{GradRegStep}
\ba{rcl}
\xx^{+} & \Def & \xx - \bigl( \nabla^2 f(\zz) + \lambda \mat{B} \bigr)^{-1}\nabla f(\xx), \qquad \lambda > 0.
\ea
\eeq
We choose parameter $\lambda$  as in Algorithm~\ref{alg:LazyGradRegNewton}:
\beq \label{LambdaChoice}
\boxed{
	\ba{rcl}
	\lambda & := & \sqrt{M \| \nabla f(\xx) \|_{*}}
	\ea
}
\eeq
where $M > 0$ is fixed.
The Newton step with gradient regularization 
can be seen as an approximation of the Cubic Newton step~\eqref{CubicStep},
which utilizes convexity of the objective. Indeed, we have by convexity:
\beq \label{GradRegRBound}
\ba{rcl}
\| \xx^{+} - \xx \| & = & \| \bigl( \nabla^2 f(\zz) + \lambda \mat{B} \bigr)^{-1}\nabla f(\xx) \|
\;\; \overset{\eqref{Convexity}}{\leq} \;\;
\frac{1}{\lambda} \| \nabla f(\xx) \|_{*},
\ea
\eeq
and so we can ensure our main bound:
\beq \label{MRBound}
\ba{rcl}
M\|\xx^{+} - \xx\| & \overset{\eqref{GradRegRBound}}{\leq} &
\frac{M}{\lambda}\|\nabla f(\xx) \|_{*}
\;\; \overset{\eqref{LambdaChoice}}{=} \;\;
\lambda,
\ea
\eeq
which justifies the actual choice of $\lambda$.

Now, we need to have analogues of Lemma~\ref{LemmaCubicGrad} and Lemma~\ref{LemmaCubicFunc} 
that we established for the cubically regularized lazy Newton step. We can prove the following.

\BL \label{LemmaGRGrad}
For any $M > 0$, it holds
\beq\label{GRGradBasic}
\ba{rcl}
\| \nabla f(\xx^{+}) \|_{*}
& \leq & 
\bigl( \frac{L}{2}\|\xx^{+} - \xx \| + \lambda + L\|\zz - \xx\| \bigr) \cdot \|\xx^{+} - \xx \|.
\ea
\eeq
Consequently, for $M \geq 3L$, we get
\beq \label{GRGrad}
\ba{rcl}
\| \nabla f(\xx^{+}) \|_{*}
 & \overset{\eqref{GRGradBasic}, \eqref{MRBound}}{\leq} & 
\bigl( \frac{7}{6} \lambda  + L \|\zz - \xx\| \bigr) \cdot \|\xx^{+} - \xx\|. 
\ea
\eeq
\EL
\proof
Indeed, we have
$$
\ba{cl}
& 
\frac{L}{2} \|\xx^{+} - \xx\|^2 \;\; \overset{\eqref{LipHessGrad}}{\geq} \;\;
\| \nabla f(\xx^{+}) - \nabla f(\xx) - \nabla^2 f(\xx)(\xx^{+} - \xx) \|_{*} \\
\\
& \geq \;\;
\| \nabla f(\xx^{+}) \|_{*} - \| \nabla f(\xx) + \nabla^2 f(\zz)(\xx^{+} - \xx) \|_{*}
- L\|\zz - \xx\| \cdot \|\xx^{+} - \xx\|,
\ea
$$
where we used triangle inequality and Lipschitz continuity of the Hessian in the last bound.
It remains to note that $\|\nabla f(\xx) + \nabla^2 f (\zz) (\xx^{+} - \xx )\|_{*} \overset{\eqref{GradRegStep}}{=} \lambda \|\xx^{+} - \xx\|$.
\qed

\BL \label{LemmaGRFunc}
For $M \geq 3L$, it holds
\beq \label{GRFunc}
\ba{rcl}
f(\xx) - f(\xx^{+}) & \geq &
\frac{\lambda}{4} \|\xx^{+} - \xx\|^2 - \frac{8L^3}{27M^2}\|\zz - \xx\|^3.
\ea
\eeq
\EL
\proof
Combining the method step with Lipschitz continuity of the Hessian, we obtain
$$
\ba{rcl}
f(\xx^{+}) & \overset{\eqref{LipHessFunc},\eqref{MRBound}}{\leq} &
f(\xx) + \la \nabla f(\xx), \xx^{+} - \xx \ra + \frac{1}{2} \la \nabla^2 f(\xx)(\xx^{+} - \xx), \xx^{+} - \xx \ra 
+ \frac{M}{18}\|\xx^{+} - \xx\|^3 \\
\\
& \overset{\eqref{GradRegStep}}{=} &
f(\xx) - \lambda \|\xx^{+} - \xx\|^2 - \la \nabla^2 f(\zz)(\xx^{+} - \xx), \xx^{+} - \xx \ra + \frac{M}{18}\|\xx^{+} - \xx\|^3 \\
\\
& & \quad + \; \frac{1}{2} \la \nabla^2 f(\xx)(\xx^{+} - \xx), \xx^{+} - \xx \ra \\
\\
& \overset{\eqref{LipHess}, \eqref{Convexity}, \eqref{MRBound}}{\leq} &
f(\xx) - \lambda \|\xx^{+} - \xx\|^2  + \frac{\lambda}{2}\|\xx^{+} - \xx\|^2 + \frac{L}{2}\| \zz - \xx\| \cdot \|\xx^{+} - \xx\|^2.
\ea
$$
Applying Young's inequality for the last term, 
$$
\ba{cl}
& \frac{L}{2}\|\zz - \xx\| \cdot \|\xx^{+} - \xx\|^2 \;\; = \;\;
\bigl(  \frac{(3M)^{2/3}}{4} \|\xx^{+} - \xx\|^2 \bigr) \cdot \bigl(  \frac{2L}{(3M)^{2/3}} \|\zz - \xx\| \bigr) \\
\\
& \leq \;\;
\frac{M}{4} \|\xx^{+} - \xx\|^3  + \frac{8 L^3}{27 M^2}\|\zz - \xx\|^3
\;\; \overset{\eqref{MRBound}}{\leq} \;\;
\frac{\lambda}{4}\|\xx^{+} - \xx\|^2 + \frac{8 L^3}{27 M^2}\|\zz - \xx\|^3.
\ea
$$
we obtain \eqref{GRFunc}.
\qed

Combining these two lemmas, we get the functional progress of one step in terms of the new gradient norm.

\BT 
Let $M \geq 3L$.
Then, for one lazy Newton step~\eqref{GradRegStep} with gradient regularization \eqref{LambdaChoice}, it holds:
\beq \label{OneStepGR}
\ba{rcl}
f(\xx) - f(\xx^{+}) & \geq & 
\frac{9}{244 \lambda} \| \nabla f(\xx^{+}) \|_{*}^2
+ \frac{M}{8} \|\xx^{+} - \xx\|^3 - \frac{3 L^3}{M^2} \|\zz - \xx\|^3.
\ea
\eeq
\ET
\proof
Using convexity of the function $t \mapsto t^2$ for $t \geq 0$ and Young's inequality, we get
\beq \label{NewGradGRTh}
\ba{cl}
& \| \nabla f(\xx^{+}) \|_{*}^2 \;\; \overset{\eqref{GRGrad}}{\leq} \;\;
\Bigl( \frac{7}{6}\lambda \|\xx^{+} - \xx\| + L \|\zz - \xx\| \cdot \|\xx^{+} - \xx\|  \Bigr)^2 \\
\\
& \; \leq \;\;
\frac{49}{18} \lambda^2 \|\xx^{+} - \xx\|^2
+ 2L^2\|\zz - \xx\|^2 \cdot \|\xx^{+} - \xx\|^2 \\
\\
& \; \leq \;\;
\frac{49}{18} \lambda^2 \|\xx^{+} - \xx\|^2
+ \frac{4 \lambda L^3}{3M^2}\|\zz - \xx\|^3
+ \frac{2M^4 \|\xx^{+} - \xx\|^6}{3 \lambda^2} \\
\\
& \overset{\eqref{MRBound}}{\leq} \;\,
\frac{61}{18} \lambda^2 \|\xx^{+} - \xx\|^2
+ \frac{4 \lambda L^3}{3M^2}\|\zz - \xx\|^3.
\ea
\eeq
Thus,
$$
\ba{rcl}
\frac{\lambda}{4}\|\xx^{+} - \xx\|^2
& \overset{\eqref{MRBound}}{\geq} &
\frac{\lambda}{8}\|\xx^{+} - \xx\|^2 + \frac{M}{8}\|\xx^{+} - \xx\|^3 \\
\\
& \overset{\eqref{NewGradGRTh}}{\geq} &
\frac{9}{244 \lambda} \| \nabla f(\xx^{+}) \|_{*}^2
+ \frac{M}{8}\|\xx^{+} - \xx\|^3
- \frac{3L^3}{61M^2} \|\zz - \xx\|^3.
\ea
$$
Substituting this bound into~\eqref{GRFunc} completes the proof.
\qed

Let us consider one phase of the algorithm. Telescoping bound~\eqref{OneStepGR}
for the first $m$ iterations and using triangle inequality,
$$
\ba{rcl}
\|\xx_0 - \xx_k\| & \leq & \sum\limits_{i = 1}^k \|\xx_i - \xx_{i - 1}\|
\;\; =: \;\; 
\sum\limits_{i = 1}^k r_i,
\ea
$$
we obtain
\BC
Let $M \geq 3mL$. Then
\beq \label{GROnePhaseBound}
\ba{rcl}
f(\xx_0) - f(\xx_m) & \geq & 
\frac{9}{244\sqrt{M}} 
\sum\limits_{k = 1}^{m} \frac{\| \nabla f(\xx_k) \|_{*}^2}{\| \nabla f(\xx_{k - 1}) \|_*^{1/2}}
+ \frac{M}{8}\sum\limits_{k = 1}^{m} r_k^3 - \frac{3L^3}{M^2} \sum\limits_{k = 1}^{m - 1}
\Bigl( \sum\limits_{i = 1}^k r_i  \Bigr)^3 \\
\\
& \overset{\eqref{TechBound}}{\geq} &
\frac{9}{244\sqrt{M}} 
\sum\limits_{k = 1}^{m} \frac{\| \nabla f(\xx_k) \|_{*}^2}{\| \nabla f(\xx_{k - 1}) \|_*^{1/2}}.
\ea
\eeq
\EC

We are ready to prove the global complexity bound for Algorithm~\ref{alg:LazyGradRegNewton}.
Let us consider the constant choice of the regularization parameter:

\beq \label{App:MGradReg}
\boxed{
\ba{rcl}
M & := & 3mL
\ea
}
\eeq

\BT \label{ThGradRegGlobalApp}
Let the objective function be convex and bounded from below: $\inf\limits_{\xx} f(\xx) \geq f^\star$.
Let the regularization parameter $M$ be fixed as in \eqref{MGradReg}.
Assume that the gradients for previous iterations $\{ \xx_i \}_{i = 1}^k$
are higher than a desired error level $\varepsilon > 0$:
\beq \label{AlgGRBigGrads}
\ba{rcl}
\| \nabla f(\xx_i) \|_{*} & \geq & \varepsilon.
\ea
\eeq
Then, the number of iterations of Algorithm~\ref{alg:LazyGradRegNewton} 
to reach  accuracy $ \| \nabla f(\xx_{k + 1}) \|_{*} \leq \varepsilon$ is at most
\beq \label{AlgGRIters}
\ba{rcl}
k & \leq & 
\cO\Bigl( \frac{\sqrt{mL} (f(\xx_0) - f^\star)}{\varepsilon^{3/2}} 
+ \ln \frac{\| \nabla f(\xx_0) \|_{*}}{\varepsilon} \Bigr).
\ea
\eeq
The total number of Hessian updates $t$ during these iterations is bounded as
\beq \label{AlgGRHess}
\ba{rcl}
t & \leq & 
\cO\Bigl( \frac{\sqrt{L} (f(\xx_0) - f^\star)}{\varepsilon^{3/2} \sqrt{m}} 
+ \frac{1}{m} \ln \frac{\| \nabla f(\xx_0) \|_{*}}{\varepsilon} \Bigr).
\ea
\eeq
\ET
\proof
Without loss of generality we assume that $k$ is a multiple of $m$:
$k = tm$. Then, for $i$th ($1 \leq i \leq t$) phase of the method, we have
\beq \label{AlgGRPhase}
\ba{rcl}
f(\xx_{m(i - 1)}) - f(\xx_{mi}) & \overset{\eqref{GROnePhaseBound}}{\geq} &
\frac{9}{244 \sqrt{M}} \sum\limits_{ m(i - 1) + 1 \, \leq \, j \, \leq \, mi }
\frac{\| \nabla f(\xx_j) \|_{*}^2}{ \| \nabla f(\xx_{j - 1}) \|_{*}^{1/2}} \\
\\
& \overset{\eqref{AlgGRBigGrads}}{\geq} &
\frac{9 \varepsilon^{3/2}}{244 \sqrt{M}} \sum\limits_{ m(i - 1) + 1 \, \leq \, j \, \leq \, mi }
\bigl(  \frac{\| \nabla f(\xx_j) \|_{*}}{ \| \nabla f(\xx_{j - 1}) \|_{*}  } \bigr)^{\frac{1}{2}}.
\ea
\eeq
Telescoping this bound for all phases and using 
inequality between arithmetic and geometric means, we get
$$
\ba{rcl}
f(\xx_0) - f^\star & \geq & f(\xx_0) - f(\xx_k)
\;\; \overset{\eqref{AlgGRPhase}}{\geq} \;\;
\frac{9 \varepsilon^{3/2}}{244 \sqrt{M}} \sum\limits_{j = 1}^k
\bigl(  \frac{\| \nabla f(\xx_j) \|_{*}}{ \| \nabla f(\xx_{j - 1}) \|_{*}  } \bigr)^{\frac{1}{2}} \\
\\
& \geq & 
\frac{9 \varepsilon^{3/2} k}{244\sqrt{M}}
\Bigl( \prod\limits_{j = 1}^k 
\frac{\| \nabla f(\xx_j) \|_{*}}{\| \nabla f(\xx_{j - 1}) \|_{*}}  \Bigr)^{\frac{1}{2k}}
\;\; = \;\;
\frac{9 \varepsilon^{3/2} k}{244\sqrt{M}}
\Bigl(  \frac{\| \nabla f(\xx_k) \|_{*} }{ \| \nabla f(\xx_0) \|_{*}  } \Bigr)^{\frac{1}{2k}} \\
\\
& \overset{\eqref{AlgGRBigGrads}}{\geq}  &
\frac{9 \varepsilon^{3/2} k}{244\sqrt{M}}
\Bigl(  \frac{  \varepsilon  }{ \| \nabla f(\xx_0) \|_{*}  } \Bigr)^{\frac{1}{2k}} 
\;\; = \;\;
\frac{9 \varepsilon^{3/2} k}{244\sqrt{M}}
\exp\Bigl(  \frac{1}{2k} \ln \frac{\varepsilon}{\| \nabla f(\xx_0) \|_{*}} \Bigr) \\
\\
& \geq & 
\frac{9 \varepsilon^{3/2} k}{244\sqrt{M}} 
\Bigl( 1 - \frac{1}{2k} \ln \frac{ \| \nabla f(\xx_0) \|_{*}}{\varepsilon} \Bigr).
\ea
$$
Rearranging the terms gives \eqref{AlgGRIters}.
Inequality \eqref{AlgGRHess} follows immediately from $t = k / m$.
\qed

We see that the global complexity bounds
for Algorithm~\ref{alg:LazyGradRegNewton}
are the same (up to an additive logarithmic term) as those ones
for Algorithm~\ref{alg:LazyCubicNewton}.
However, each iteration of Algorithm~\ref{alg:LazyGradRegNewton}
is much easier to implement since it involves just one standard matrix inversion.
\section{Proofs for Section~\ref{SectionLocal}}
 We assume here that our objective is strongly convex
 with some parameter $\mu > 0$.
\beq \label{SrongConvex}
\ba{rcl}
\nabla^2 f(\xx) & \succeq & \mu \mat{B}, \qquad \forall \xx \in \R^d.
\ea
\eeq

By strong convexity~\eqref{SrongConvex}, we have a bound for the distance to the optimum
$\xx^\star \Def \argmin_{\xx} f(\xx)$ using the gradient norm (see, e.g. \cite{nesterov2018lectures}),
for all $\xx \in \R^d$:
\beq \label{DistanceStrongly}
\ba{rcl}
\|\xx - \xx^\star\| & \leq &
\frac{1}{\mu} \| \nabla f(\xx) \|_{*}.
\ea
\eeq

Now, let us look at one lazy Cubic Newton step $\mat{T} = \mat{T}_{M}(\xx, \zz)$ with some $M \geq 0$,
for which we have:
\beq \label{ROldGradBound}
\ba{rcl}
r & \equiv & \|\mat{T} - \xx\|
\;\; \overset{\eqref{CubicStep}}{=} \;\;
\| \bigl(  \nabla^2 f(\zz) + \frac{Mr}{2}\mat{B} \bigr)^{-1} \nabla f(\xx) \|_{*}
\;\; \overset{\eqref{SrongConvex}}{\leq} \;\;
\frac{1}{\mu} \| \nabla f(\xx) \|_{*}.
\ea
\eeq
Note that $M = 0$ corresponds to the classical pure Newton step
with a lazy Hessian.

Then, applying Lemma~\ref{LemmaCubicGrad}
with triangle inequality, we get
\beq \label{SuperlinearNewGrad}
\ba{rcl}
\| \nabla f(\mat{T}) \|_{*} & \overset{\eqref{NewGradBound}}{\leq} &
\frac{M + L}{2} r^2 + Lr \|\zz - \xx\| \;\;
\leq \;\;
\frac{M + L}{2} r^2 + Lr \|\xx - \xx^\star \| + Lr \|\zz - \xx^\star\| \\
\\
& \overset{\eqref{DistanceStrongly}, \eqref{ROldGradBound}}{\leq} &
\frac{M + 3L}{2\mu^2}
\| \nabla f(\xx) \|_{*}^2
+ \frac{L}{\mu^2} \| \nabla f(\xx) \|_{*} \cdot \| \nabla f(\zz) \|_{*}.
\ea
\eeq
This inequality leads to a superlinear convergence of the method.

\BT \label{ThLocalApp}
Let $M \geq 0$.
Assume that initial gradient is small enough:
\beq \label{InitGradSuperlinear}
\ba{rcl}
\| \nabla f(\xx_0) \|_{*} & \leq & \frac{\mu^2}{2(M + 3L)}.
\ea
\eeq
Then, for the iterations of Algorithm~\ref{alg:LazyCubicNewton}, we have
a superlinear convergence for the gradient norms:
\beq \label{CubicSuperlinear}
\ba{rcl}
\| \nabla f(\xx_k) \|_{*}
& \leq & 
\frac{\mu^2}{M + 3L}
\bigl( \frac{1}{2}  \bigr)^{(1 + m)^{\pi(k)}(1 + k \Mymod m) },
\qquad k \geq 0,
\ea
\eeq
where $\pi(k)$ is defined by \eqref{PiDef}.
\ET
\proof
According to~\eqref{SuperlinearNewGrad}, for one iteration 
of Algorithm~\ref{alg:LazyCubicNewton}, it holds
$$
\ba{rcl}
\| \nabla f(\xx_{k + 1}) \|_{*}
& \leq & \frac{M + 3L}{2\mu^2} \|\nabla f(\xx_k)\|_{*}^2
+ \frac{L}{\mu^2} \| \nabla f(\xx_k) \|_{*} \cdot \| \nabla f(\xx_{\pi(k)}) \|_{*}.
\ea
$$
Let us multiply this inequality by $c := \frac{M + 3L}{2\mu^2}$
and denote $s_k \Def c \| \nabla f(\xx_k) \|_{*}$. This yields
the following recursion
\beq \label{LocalMainRec}
\ba{rcl}
s_{k + 1} & \leq & s_k^2 + s_k s_{\pi(k)},
\ea
\eeq
and by initial condition \eqref{InitGradSuperlinear}, we have $s_0 \leq \frac{1}{4}$.

Required inequality from our claim is
\beq \label{LocalRecSuff}
\ba{rcl}
s_{tm + i} & \leq & \bigl( \frac{1}{2} \bigr)^{(1 + i)(1 + m)^t + 1},
\ea
\eeq
for any $t \geq 0$ (phases of the method) and $0 \leq i \leq m$.
We prove it by induction. Indeed,
$$
\ba{rcl}
s_{tm + i + 1}
& \overset{\eqref{LocalMainRec}}{\leq} &
s_{tm + i}^2 + s_{tm + i} s_{tm}
\;\; \overset{\eqref{LocalRecSuff}}{\leq} \;\;
\bigl( \frac{1}{2} \bigr)^{2(1 + i)(1 + m)^t + 2}
+ 
\bigl( \frac{1}{2} \bigr)^{(1 + i)(1 + m)^t + 1 + (1 + m)^t + 1} \\
\\
& \leq & 
\bigl( \frac{1}{2} \bigr)^{(2 + i)(1 + m)^t + 2}
+
\bigl( \frac{1}{2} \bigr)^{(2 + i)(1 + m)^t + 2}
\;\; = \;\;
\bigl( \frac{1}{2} \bigr)^{(2 + i)(1 + m)^t + 1},
\ea
$$
which is \eqref{LocalRecSuff} for the next index.
\qed

\BC
Combining both Theorem~\ref{ThGlobCompl} and \ref{ThLocal}, we conclude
that for minimizing a strongly convex function
by Algorithm~\ref{alg:LazyCubicNewton} with regularization parameter 
$M$ given by \eqref{MCubicNewton} and starting from an arbitrary $\xx_0$,
we need to do
$$
\ba{rcl}
k & \leq & 
\cO\Bigl( 
\frac{m^2 L^2 (f(\xx_0) - f^\star)}{\mu^3}
\, + \, \frac{1}{\ln(1 + m)} \ln \ln \frac{\mu^2}{mL\varepsilon}
\Bigr)
\ea
$$
lazy steps to achieve $\| \nabla f(\xx_k) \|_{*} \leq \varepsilon$.
\EC

Let us also analyze the local behaviour of the lazy Newton steps with gradient regularization,
that we perform in Algorithm~\ref{alg:LazyGradRegNewton}.
Note that by strong convexity~\eqref{SrongConvex} we have 
the following bound for the length of one step, which is
\beq \label{GROldGradBound}
\ba{rcl}
\| \xx^{+} - \xx \| & \overset{\eqref{GradRegStep}}{=} &
\| \bigl(  \nabla^2 f(\zz) + \lambda \mat{B}  \bigr)^{-1} \nabla f(\xx) \|
\;\; \overset{\eqref{SrongConvex}}{\leq} \;\; \frac{1}{\mu} \| \nabla f(\xx) \|_{*}.
\ea
\eeq
Therefore, we can estimate the norm of the gradient at new point
taking into account the actual choice of the regularization parameter~\eqref{LambdaChoice}.
Our reasoning works for any $M \geq 0$.

Applying Lemma~\ref{LemmaGRGrad} and triangle inequality, we get
\beq \label{GRNewGradLocal}
\ba{cl}
& \| \nabla f(\xx^{+}) \|_{*} \;\; \overset{\eqref{GRGradBasic}}{\leq} \;\;
\bigl(  \frac{L}{2} \|\xx^{+} - \xx\| + \lambda + L\|\zz - \xx\| \bigr) \cdot \|\xx^{+} - \xx\| \\
\\
& \;\;\,  \overset{\eqref{GROldGradBound}}{\leq} \;\;\, \;
\bigl(  \frac{L}{2\mu} \| \nabla f(\xx) \|_{*} + \lambda + L\|\zz - \xx\| \bigr) \cdot 
\frac{1}{\mu} \| \nabla f(\xx) \|_{*} \\
\\
& \;\;\, \; \leq \;\; \;\;\, 
\bigl(  \frac{L}{2\mu} \| \nabla f(\xx) \|_{*} + \lambda + L\|\xx - \xx^\star\|
+ L\|\zz - \xx^\star\| \bigr) \cdot 
\frac{1}{\mu} \| \nabla f(\xx) \|_{*}  \\
\\
& \overset{\eqref{DistanceStrongly},\eqref{LambdaChoice}}{\leq} \;
\frac{3L}{2\mu^2} \| \nabla f(\xx) \|_{*}^2
+ \frac{\sqrt{M}}{\mu} \| \nabla f(\xx) \|_{*}^{3/2}
+ \frac{L}{\mu^2} \| \nabla f(\xx) \|_{*} \cdot \| \nabla f(\zz) \|_{*} 
\ea
\eeq

We see that the power of the gradient norm for the second term in the right hand side is $3 / 2$,
which is slightly worse than $2$, that is in the cubic regularization (compare with \eqref{SuperlinearNewGrad}).
However, it is still a local \textit{superlinear} convergence, which is extremely fast from the practical perspective.

\BT \label{ThGRLocalApp}
Let $M \geq 0$.
Assume that initial gradient is small enough:
\beq \label{GRInitGradSuperlinear}
\ba{rcl}
\| \nabla f(\xx_0) \|_{*} & \leq & 
\frac{\mu^2}{2^5( 3L + 4M)}.
\ea
\eeq
Then, for the iterations of Algorithm~\ref{alg:LazyGradRegNewton},
we have a superlinear convergence for the gradient norms:
\beq \label{GRSuperlinear}
\ba{rcl}
\| \nabla f(\xx_k) \|_{*}
& \leq & 
\frac{\mu^2}{2^3 (3L + 4M)}
\bigl( \frac{1}{2} \bigr)^{2( 1 \, + \, m / 2 )^{\pi(k)} (1 + (k \Mymod m) / 2) },
\qquad k \geq 0,
\ea
\eeq
where $\pi(k)$ is defined by \eqref{PiDef}.
\ET
\proof
According to \eqref{GRNewGradLocal}, for one iteration of Algorithm~\ref{alg:LazyGradRegNewton},
it holds
$$
\ba{rcl}
\| \nabla f(\xx_{k + 1}) \|_{*} & \leq &
\frac{3L}{\mu^2} \| \nabla f(\xx_k) \|_{*}^2
+ \frac{\sqrt{M}}{\mu} \| \nabla f(\xx_k) \|_{*}^{3/2}
+ \frac{L}{\mu^2} \| \nabla f(\xx_k) \|_{*} \cdot \| \nabla f(\xx_{\pi(k)}) \|_{*}.
\ea
$$
Let us multiply this inequality by $c := \frac{2(4M + 3L)}{\mu^2}$
and denote $s_k \Def c \| \nabla f(\xx_k) \|_{*}$. 
This yields the following recursion (compare with \eqref{LocalMainRec}):

\beq \label{GRLocalMainRec}
\ba{rcl}
s_{k + 1} & \leq & \frac{1}{2}\bigl( s_k^2 +  s_k^{3/2}\bigr) + s_k s_{\pi(k)},
\ea
\eeq
and by initial condition~\eqref{GRInitGradSuperlinear}, we have $s_0 \leq \frac{1}{2^4}$.

Let us prove by induction that
\beq \label{GR_local_stm}
\ba{rcl}
s_{tm} & \leq & \bigl(  \frac{1}{2} \bigr)^{2( 1 + (1 + m/2)^t )},
\ea
\eeq
for any $t \geq 0$ (phases of the method), and
\beq \label{GR_local_stmi}
\ba{rcl}
s_{tm + i} & \leq & 
\bigl( \frac{1}{2} \bigr)^{i(1 + m / 2)^t} s_{tm},
\ea
\eeq
for $0 \leq i \leq m$. Then \eqref{GRSuperlinear} clearly follows.

Note that from \eqref{GR_local_stm} we have
\beq \label{GR_local_sqrt}
\ba{rcl}
\frac{1}{2}s_{tm}^{1/2} +  \frac{3}{2} s_{tm} & \leq &
2 s_{tm}^{1/2}
\;\; \overset{\eqref{GR_local_stm}}{\leq} \;\;
\bigl( \frac{1}{2} \bigr)^{(1 + m/2)^t}.
\ea
\eeq
Then, assuming that \eqref{GR_local_stm} and \eqref{GR_local_sqrt}
holds for the current iterate, we get
$$
\ba{rcl}
s_{tm + i + 1}
& \overset{\eqref{GRLocalMainRec}}{\leq} &
\frac{1}{2}
\bigl( 
s_{tm + i}^2 +
s_{tm + i}^{3/2}
\bigr) + s_{tm + i} s_{tm}
\;\; \overset{\eqref{GR_local_stmi}, \eqref{GR_local_stm}}{\leq} \;\;
\bigl( \frac{1}{2} \bigr)^{i(1 + m/2)^t} s_{tm}
\bigl(  \frac{1}{2} s_{tm}^{1/2} + \frac{3}{2} s_{tm} \bigr) \\
\\
& \overset{\eqref{GR_local_sqrt}}{\leq} &
\bigl( \frac{1}{2} \bigr)^{(i + 1)(1 + m/2)^t} s_{tm},
\ea
$$
which is \eqref{GR_local_stmi} for the next index.

It remains to observe that substituting $i := m$ into \eqref{GR_local_stmi}
we obtain \eqref{GR_local_stm} for the next phase. Indeed,
$$
\ba{rcl}
s_{(t + 1)m}
& = & 
s_{tm + m}
\;\; \overset{\eqref{GR_local_stmi}}{\leq} \;\;
\bigl( \frac{1}{2} \bigr)^{m(1 + m/2)^t} s_{tm} \\
\\
& \overset{\eqref{GR_local_stm}}{\leq} &
\bigl( \frac{1}{2} \bigr)^{m(1 + m/2)^t + 2 + 2(1 + m/2)^t} 
\;\; = \;\;
\bigl( \frac{1}{2} \bigr)^{2 + 2(1 + m/2)^{t + 1}}. \QF
\ea
$$

\BC
Combining both Theorem~\ref{ThGradRegGlobal} and \ref{ThGRLocal}, we conclude
that for minimizing a strongly convex function
by Algorithm~\ref{alg:LazyGradRegNewton} with regularization parameter 
$M$ given by \eqref{MGradReg} and starting from an arbitrary $\xx_0$,
we need to do
$$
\ba{rcl}
k & \leq & 
\cO\Bigl( 
\frac{m^2 L^2 (f(\xx_0) - f^\star)}{\mu^3}
\, + \, \ln \frac{mL \|\nabla f(\xx_0) \|_{*}}{\mu^2}
\, + \, \frac{1}{\ln(1 + m /2)} \ln \ln \frac{\mu^2}{mL\varepsilon}
\Bigr)
\ea
$$
lazy steps to achieve $\| \nabla f(\xx_k) \|_{*} \leq \varepsilon$.
\EC
\section{Proofs for Section~\ref{SectionImplementation} and the Adaptive Regularized Newton}
\subsection{Proof of the General Case}
\label{GC}
The key to our analysis of Algorithm~\ref{alg:AdaptiveCubicNewton} is inequality \eqref{Telescoped} on $m$ consecutive steps
of the method. Note that it holds for \textit{any} value of $M$ that is sufficiently big.
Hence, there is no necessity to fix the regularization parameter. We are going to 
change its value adaptively while checking the condition on the functional progress
after $m$ steps of the method.

The value $M_0$ is just an \textit{initial guess} for the regularization constant,
which can be further both increased or decreased dynamically. 
This process is well defined. Indeed, due to \eqref{Telescoped}, the stopping condition
is satisfied as long as $M_t \geq 2^{8} 3^5 mL$. Hence, we ensure that 
all regularization coefficients always satisfy the following bound:
\beq \label{AdAlgMBound}
\ba{rcl}
M_t & \leq & \max\{ 2M_0, 2^9 3^5 mL \}.
\ea
\eeq
Substituting this bound into the stopping condition, and telescoping it for the first $t \geq 0$ phases, we get
\beq \label{AdAlgGrads}
\ba{rcl}
f(\xx_0) - f^\star & \geq & 
f(\xx_0) - f(\xx_{tm}) \;\; \geq \;\;
\frac{1}{ \sqrt{\max\{ 2M_0, 2^9 3^5 mL \}} }
\sum\limits_{i = 1}^{tm} \| \nabla f(\xx_i) \|_{*}^{3/2}.
\ea
\eeq
Thus, we obtain the following global complexity guarantee.

\BT \label{ThAdaptiveCubicGloblalApp}
Let the objective function be bounded from below: $\inf\limits_{\xx} f(\xx) \geq f^\star$ and
let the gradients for previous iterates $\{ \xx_i \}_{i = 1}^{tm}$ 
during $t$ phases of Algorithm~\ref{alg:AdaptiveCubicNewton}
are higher than a desired error level $\varepsilon > 0$:
\beq \label{AdAlgGLB}
\ba{rcl}
\| \nabla f(\xx_i) \|_{*} & \geq & \varepsilon.
\ea
\eeq
Then, the number of phases of Algorithm~\ref{alg:AdaptiveCubicNewton}
to reach accuracy $\| \nabla f(\xx_{tm + 1}) \|_{*} \leq \varepsilon$
is at most
\beq \label{AdAlgPhases}
\ba{rcl}
t & \leq & \cO\Bigl( 
\sqrt{
\max\bigl\{ \frac{M_0}{m}, L \bigr\} } \cdot
 \frac{f(\xx_0) - f^\star}{\varepsilon^{3/2}\sqrt{m}}  \Bigr).
\ea
\eeq
The total number $N$ of gradient calls during these phases is bounded as
\beq \label{App:AdAlgTotal}
\ba{rcl}
N & \leq & \bigl(2 t + \log_2 \frac{M_t}{M_0} \bigr)m
\;\; \overset{\eqref{AdAlgMBound}}{\leq} \;\;
2tm
\, + \,
 \max\{1, \log_2 \frac{2^9 3^5 mL}{M_0}  \} m.
\ea
\eeq
\ET 
\proof
Bound~\eqref{AdAlgPhases} follows directly from \eqref{AdAlgGrads}
and \eqref{AdAlgGLB}.

To prove \eqref{App:AdAlgTotal}, we denote by $n_i$ the number
of times the \textbf{do \ldots until} loop is performed at phase $0 \leq i \leq {t - 1}$.
By our updates, it holds
\beq \label{AdAlgPowers}
\ba{rcl}
2^{n_i - 2} = \frac{M_{i + 1}}{M_i}.
\ea
\eeq
It remains to note that
$$
\ba{rcl}
N & = & \biggl(   \sum\limits_{i = 0}^{t - 1} n_i \biggr) m
\;\; \overset{\eqref{AdAlgPowers}}{=} \;\;
\biggl(  2t + \sum\limits_{i = 0}^{t - 1} \log_2 \frac{M_{i + 1}}{M_i} \biggr) m
\;\; = \;\;
\biggl(  2t + \log_2 \frac{M_t}{M_0} \biggr) m. \QF
\ea
$$

\BR
According to Theorem \ref{ThAdaptiveCubicGloblal}, 
the global complexity of the adaptive algorithm
is the same as that one given by Theorem~\ref{ThGlobCompl} 
for the method with a fixed regularization constant.
Due to \eqref{App:AdAlgTotal}, the average number
of tries of different $M_t$ per one phase is only \underline{two}.
\ER

\subsection{Adaptive Search for the Regularized Newton}
\label{ASCV}

Let us present an adaptive version of Algorithm~\ref{alg:LazyGradRegNewton},
which is suitable for solving convex minimization problems (Algorithm~\ref{alg:AdaptiveGradRegNewton} below).

\begin{algorithm}[h!]
	\caption{Adaptive Regularized Newton with Lazy Hessians}
	\label{alg:AdaptiveGradRegNewton}
	\begin{algorithmic}[1]
		\STATE {\bfseries Input:} $\xx_0 \in \R^d$, $m \geq 1$. Fix some $M_0 > 0$.
		\FOR{$t=0,1,\dotsc$}
		\STATE Compute snapshot Hessian $\nabla^2 f( \xx_{tm} )$
		\REPEAT
		\STATE Update $M_t = 2 \cdot M_t$
		\FOR{$i=1,\ldots,m$}
		\STATE Denote $k = tm + i - 1$ and set  $\lambda_k = \sqrt{M_t \| \nabla f(\xx_k) \|_{*}}$
		\STATE Compute lazy Newton step 
		$\xx_{k + 1} = \bigl( \nabla^2 f(\xx_{tm}) + \lambda_k \mat{B}\bigr)^{-1} \nabla f(\xx_{k})$
		\ENDFOR
		\UNTIL {$f(\xx_{tm}) - f(\xx_{tm + m}) \; \geq \; 
			\sum_{i = 1}^m  \frac{1}{\lambda_{tm + i - 1}} \| \nabla f(\xx_{tm + i})   \|_{*}^{2} $}
		\STATE Set $M_{t + 1} = \frac{1}{4} \cdot M_t$
		\ENDFOR
	\end{algorithmic}
\end{algorithm}

 Repeating the same reasoning as in section~\ref{GC} , we straightforwardly obtain
 the same global guarantees for this process as in Theorem~\ref{ThGradRegGlobal}
 for the method with a fixed Lipschitz constant.
The cost of adaptive search again is only one extra try of regularization parameter
in average per phase.

\newpage

\section{Composite Convex Optimization}
\label{SectionCCase}

Let us consider \textit{composite} convex optimization problems
in the following form:
\beq \label{CompositeProblem}
\ba{rcl}
\min\limits_{\xx \in \dom \psi} \Bigl\{ 
F(\xx) & \Def & f(\xx) + \psi(\xx)
\Bigr\},
\ea
\eeq
where $f: \R^{d} \to \R$ is a several times differentiable convex function, 
with Lipschitz continuous Hessian~\eqref{OldLipHess},
and $\psi: \R^{d} \to \R \cup \{ +\infty \}$
is a proper closed convex function, that can be non-differentiable. However, we assume 
that it has a \textit{simple} structure
which allows to solve efficiently corresponding minimization subproblems
that involve $\psi(\cdot)$. 

In this section, we demonstrate how to generalize our analysis
of the lazy Newton steps onto this class of problems.
For example, $\psi$ can be an indicator of 
a simple closed convex set $Q \subset \R^d$:
$$
\ba{rcl}
\psi(\xx) & := & \begin{cases}
0, \quad & \xx \in Q, \\
+\infty, \quad & \text{otherwise},
\end{cases}
\ea
$$
in which case problem \eqref{CompositeProblem} becomes 
\textit{constrained} optimization problem:
$$
\ba{c}
\min\limits_{\xx \in Q} f(\xx).
\ea
$$

\subsection{Lazy Cubic Newton for Composite Problems}

We replace the lazy Cubic Newton step \eqref{StepDef}
by the following one:
\beq \label{CCstep}
\ba{rcl}
\mat{T}_{M}(\xx, \zz) & = & \argmin\limits_{\yy \in \dom \psi}
\Bigl\{ 
\la \nabla f(\xx), \yy - \xx \ra + \frac{1}{2} \la \nabla^2 f(\zz)(\yy - \xx), \yy - \xx \ra
+ \frac{M}{6}\|\yy - \xx\|^3 + \psi(\yy)
\Bigr\}.
\ea
\eeq
Since we assume here that both $f$ and $\psi$ \textit{convex},
the solution to \eqref{CCstep} always exists and is unique due to the uniform
convexity of the cubic regularizer \cite{nesterov2008accelerating}.
Then, we can use this point exactly the same way as 
the basic step~\eqref{StepDef}
in Algorithm~\ref{alg:LazyCubicNewton} and Algorithm~\ref{alg:AdaptiveCubicNewton}.

Let us present the main inequalities that are needed for its theoretical analysis.

The stationary condition for $\mat{T} := \mat{T}_M(\xx, \zz)$ is as follows, for any $\yy \in \dom \psi$,
\beq \label{StatCond}
\ba{rcl}
\la \nabla f(\xx) + \nabla^2 f(\zz) (\mat{T} - \xx) + \frac{M}{2} r \mat{B}(\mat{T} - \xx) , \yy - \mat{T} \ra
+ \psi(\yy) 
& \geq & 
\psi(\mat{T}),
\ea
\eeq
where $r := \| \mat{T} - \xx \|$ as always. In other words, we have that
\beq \label{SubgrDef}
\ba{rcl}
\psi'(\mat{T}) & \Def &
- \nabla f(\xx) - \nabla^2 f(\zz) (\mat{T} - \xx) - \frac{M}{2} r
\mat{B}(\mat{T} - \xx)
\;\; \in \;\; \partial \psi(\mat{T}).
\ea
\eeq
Correspondingly, we denote 
$$
\ba{rcl}
F'(\mat{T}) & \Def & \nabla f(\mat{T}) + \psi'(\mat{T}) \;\; \in \;\; \partial F(\mat{T}).
\ea
$$
We can work with these objects in a similar way as with the new gradient of $f$ in
the basic non-composite case. Thus,
we can ensure the following inequalities,
employing the Lipschitzness of the Hessian:
\beq \label{CompositeNewGrad}
\ba{rcl}
\frac{Lr^2}{2} 
& \overset{\eqref{LipHessGrad}}{\geq} & 
\| \nabla f(\mat{T}) - \nabla f(\xx) - \nabla^2 f(\zz)(\mat{T} - \xx) \|_{*} \\
\\
& \overset{\eqref{SubgrDef}}{=} &
\| F'(\mat{T}) + \frac{Mr}{2} \mat{B}(\mat{T} - \xx)
+ (\nabla^2 f(\zz) - \nabla^2 f(\xx)) (\mat{T} - \xx) \|_{*} \\
\\
& \geq &
\| F'(\mat{T}) \|_{*} - \frac{Mr^2}{2} - r \| \nabla^2 f(\zz) - \nabla^2 f(\xx) \| \\
\\
& \overset{\eqref{OldLipHess}}{\geq} & 
\| F'(\mat{T}) \|_{*} - \frac{Mr^2}{2}  - r L \| \zz - \xx\|.
\ea
\eeq
Hence, rearranging the terms we can prove the analogue of Lemma~\ref{LemmaCubicGrad}
for the composite case \eqref{CompositeProblem}:
\BL \label{LemmaCompositeCubicGrad}
For any $M > 0$, it holds that
\beq \label{NewGradCompositeBound}
\ba{rcl}
\| F'(\mat{T}) \|_{*} & \leq & \frac{M + L}{2} r^2 + L r\|\zz - \xx\|.
\ea
\eeq
\EL

Secondly, for the objective function value at new point, we get
$$
\ba{rcl}
F(\mat{T}) & \Def & f(\mat{T}) + \psi(\mat{T}) \\
\\
& \overset{\eqref{LipHessFunc}}{\leq} &
f(\xx) + \la \nabla f(\xx), \mat{T} - \xx \ra
+ \frac{1}{2} \la \nabla^2 f(\xx)(\mat{T} - \xx), \mat{T} - \xx \ra
+ \frac{L}{6} r^3 + \psi(\mat{T}) \\
\\
& \leq & 
F(\xx)
+ \la \nabla f(\xx), \mat{T} - \xx \ra
+ \frac{1}{2} \la \nabla^2 f(\xx)(\mat{T} - \xx), \mat{T} - \xx \ra
+ \frac{L}{6} r^3 + \la \psi'(\mat{T}), \mat{T} - \xx \ra \\
\\
& \overset{\eqref{SubgrDef}}{=} &
F(\xx)
+ \frac{1}{2} \la \nabla^2 f(\xx)(\mat{T} - \xx), \mat{T} - \xx \ra
- \la \nabla^2 f(\zz)(\mat{T} - \xx), \mat{T} - \xx \ra
- \frac{3M - L}{6} r^3 \\
\\
& \overset{(*)}{\leq} &
F(\xx) + \frac{1}{2} \la ( \nabla f(\xx) - \nabla f(\zz) ) (\mat{T} - \xx), \mat{T} - \xx \ra
- \frac{3M - L}{6} r^3 \\
\\
& \overset{\eqref{OldLipHess}}{\leq} &
F(\xx) + \frac{L}{2} r^2 \|\zz - \xx \| - \frac{3M - L}{6} r^3,
\ea
$$
where we used convexity of $f(\cdot)$ in $(*)$.
Note that the form of this inequality exactly the same as in the proof
of Lemma~\ref{LemmaCubicFunc}.
Hence, we can establish its analogue fo the composite case~\eqref{CompositeProblem}:

\BL \label{LemmaCompositeCubicFunc}
For any $M > 0$, it holds that
\beq \label{NewFuncCompositeProgress}
\ba{rcl}
F(\xx) - F(\mat{T}) & \geq & \frac{11M - 4L}{24} r^3
- \frac{32L^3}{3M^2}\|\zz - \xx\|^3.
\ea
\eeq
\EL

Therefore, having these two main lemmas established, we can prove 
the global rates for the composite Cubic Newton with Lazy Hessian updates
using similar technique that were used before.
Note that in the composite case, we ensure the convergence
in terms of the \textit{subgradient} norm: $\| F'(\xx_{k}) \|_{*} \to 0$
with $k \to \infty$.

Finally, let justify local superlinear convergence, when the smooth component
of the objective is \textit{strongly convex}: 
$\nabla^2 f(\xx) \succeq \mu \mat{B}, \forall x \in \R^d$.
Then, we have the following inequality for the whole objective:
\beq \label{CompositeSC}
\ba{rcl}
F(\yy) & \geq & F(\xx) + \la \ss_{\xx}, \yy - \xx \ra + \frac{\mu}{2}\| \yy - \xx\|^2,
\qquad \forall \xx, \yy \in \dom F, \quad  \ss_{\xx} \in \partial F(\xx).
\ea
\eeq
Substituting $\xx := \xx^{\star}$ into \eqref{CompositeSC}, we get, for any $\yy \in \dom F$:
\beq \label{FstarBound}
\ba{rcl}
F(\yy) - F^{\star} & \geq & \frac{\mu}{2} \| \yy - \xx^{\star} \|^2.
\ea
\eeq
At the same time, minimizing the left and right hand sides of \eqref{CompositeSC}
with respect to $\yy$ independently, we obtain
\beq \label{FGradBound}
\ba{rcl}
F^{\star} & \geq & F(\xx) - \frac{1}{2\mu} \| \ss_{\xx} \|_{*}^2.
\ea
\eeq
Thus, combining \eqref{FstarBound} and \eqref{FGradBound} together, we ensure
the following standard inequality, for any $\xx \in \dom F$ and $\ss_{\xx} \in \partial F(\xx)$:
\beq \label{FStrongDist}
\ba{rcl}
\frac{1}{\mu} \| \ss_{\xx} \|_{*} & \geq & \| \xx - \xx^{\star} \|.
\ea
\eeq
Further, for one lazy composite cubic step $\xx \mapsto \mat{T} := \mat{T}_M(\xx, \zz)$, we obtain
$$
\ba{rcl}
\mu r^2 & \leq & 
\la \nabla^2 f(\zz) (\mat{T} - \xx), \mat{T} - \xx \ra 
\;\; \overset{\eqref{SubgrDef}}{=} \;\;
\la \nabla f(\xx) + \psi'(\mat{T}), \xx - \mat{T} \ra
- \frac{M}{2} r^3 \\
\\
& \overset{(*)}{\leq} & \la \ss_{\xx}, \xx - \mat{T} \ra - \frac{M}{2} r^3
\;\; \leq \;\;
r \| \ss_{\xx} \|_{*},
 \ea
$$
where we used convexity of $\psi$ in $(*)$. Hence,
\beq \label{FHessLocal}
\ba{rcl}
\frac{1}{\mu} \| \ss_{\xx} \|_{*} & \geq & r.
\ea
\eeq
It remains to use Lemma~\ref{LemmaCompositeCubicGrad}, which gives,
for any $M \geq 0$:
$$
\ba{rcl}
\| F'(\mat{T}) \|_{*} & \overset{\eqref{NewGradCompositeBound}}{\leq} &
\frac{M + L}{2} r^2 + Lr \| \zz - \xx \|
\;\; \leq \;\;
\frac{M + L}{2} r^2
+ Lr \| \xx - \xx^{\star} \|
+ Lr \| \zz - \xx^{\star} \| \\
\\
& \overset{\eqref{FStrongDist}, \eqref{FHessLocal}}{\leq} & 
\frac{M + 3L}{2\mu^2} \| F'(\xx)\|_{*}^2
\; + \; \frac{L}{\mu} \| F'(\xx) \|_{*} \cdot \| F'(\zz) \|_*.
\ea
$$
That ensures local superlinear convergence
in terms of the subgradient norm
for our method in the composite case 
(compare with \eqref{SuperlinearNewGrad}).

\section{Extra Experiments}
\label{SectionEE}

In this section, we include additional evaluation
of our methods on several optimization problems:
\textit{Logistic Regression with $\ell_2$-regularization},
\textit{Logistic Regression with a non-convex regularizer},
and training a small \textit{Diagonal Neural Network} model.
In all our experiments we observe that the use of
the lazy Hessian updates significantly improve the performance 
of the Cubic Newton method
in terms of the total computational cost.
We also show the convergence of the  classic Gradient Method (GM)
as a natural baseline.
We use a constant regularization parameter $M$ (correspondingly the stepsize in the Gradient Method), that we choose for each method separately to optimize its performance.

\subsection{Convex Logistic Regression}

\begin{figure}[h!]
	\centering
	\includegraphics[width=0.42\textwidth ]{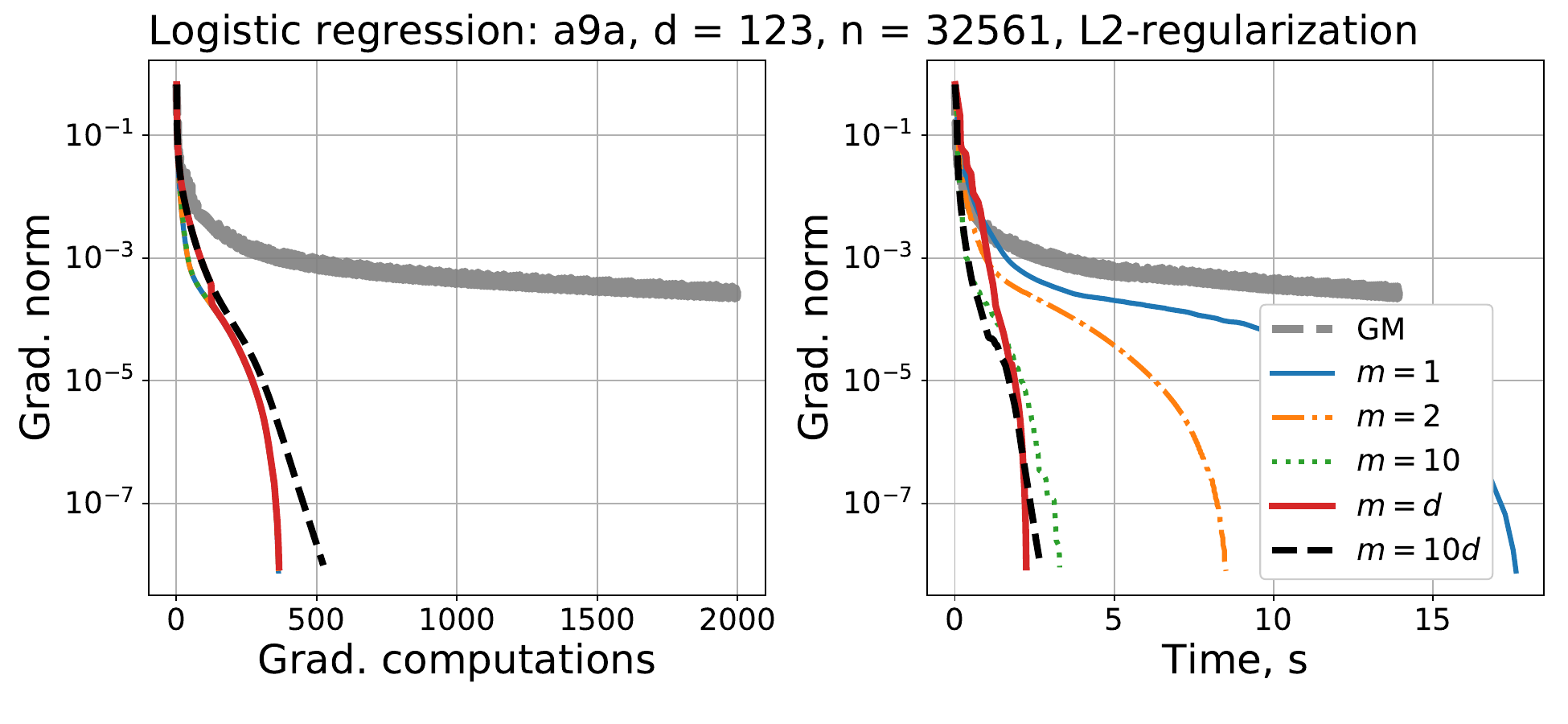} \\
	\hspace*{-5pt}
	\caption{\small Logistic Regression with $\ell_2$-regularization trained on the \texttt{a9a} dataset.  The Cubic Newton method with Lazy Hessian updates ($m = d$) shows the best overall performance.  }
	\label{fig:LogRegConvex}
\end{figure}

The objective in this problem has the following form:
$$
\ba{rcl}
f(\xx) & = & 
\frac{1}{n} \sum\limits_{i = 1}^n
\log (1 + e^{-y_i \la \aa_i, \xx \ra} ) + \frac{\lambda}{2}\| \xx\|^2,
\qquad \xx \in \R^d,
\ea
$$
where $\aa_1, \ldots, \aa_n \in \R^d$ are the feature vectors
and $y_1, \ldots, y_n \in \{ -1, 1\}$ are the labels, given by the dataset~\footnote{\href{https://www.csie.ntu.edu.tw/~cjlin/libsvmtools/datasets/}{www.csie.ntu.edu.tw/~cjlin/libsvmtools/datasets/}.}
and $\lambda > 0$ is the regularization parameter, which we fix as $\lambda = \frac{1}{n}$.
The results are shown in Figure~\ref{fig:LogRegConvex}.

\newpage
\subsection{Non-convex Logistic Regression}

\begin{figure}[h!]
	\centering
	\includegraphics[width=0.42\textwidth ]{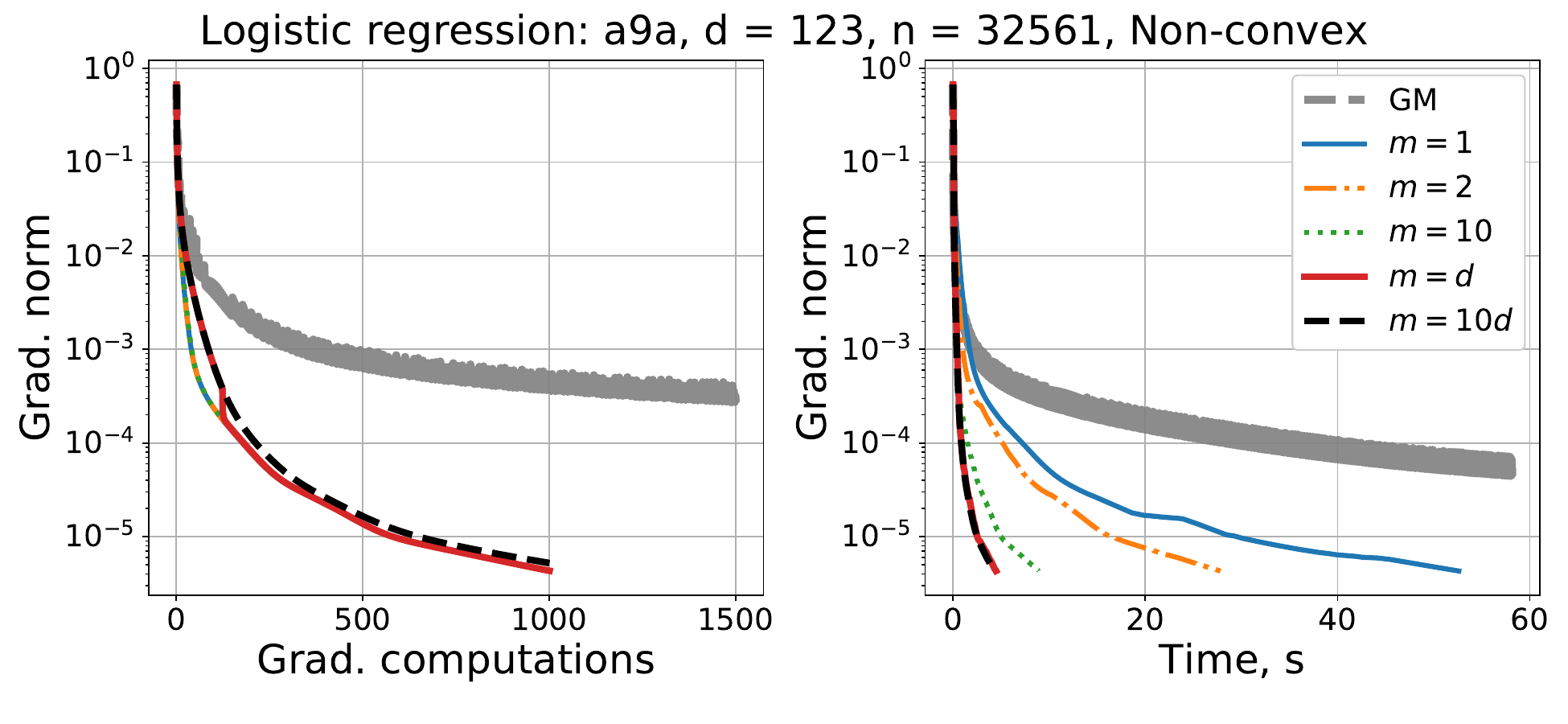} \\
	\hspace*{-5pt}
	\caption{\small Logistic Regression with non-convex regularizer trained on the \texttt{a9a} dataset. The Cubic Newton method with Lazy Hessian updates ($m = d$) shows the best overall performance.}
	\label{fig:LogRegNonConvex}
\end{figure}

In the following experiment, we consider the Logistic Regression model
with \textit{non-convex} regularizer:
$$
\ba{rcl}
f(\xx) & = & 
\frac{1}{n} \sum\limits_{i = 1}^n
\log (1 + e^{-y_i \la \aa_i, \xx \ra} ) 
+ \lambda \sum\limits_{i = 1}^d \frac{x_i^2}{1 + x_i^2},
\ea
$$
that was also studied in \cite{kohler2017sub}, and we set $\lambda = \frac{1}{n}$.
The results are shown in Figure~\ref{fig:LogRegNonConvex}.

\subsection{Non-convex Diagonal Neural Network}

\begin{figure}[h!]
	\centering
	\includegraphics[width=0.42\textwidth ]{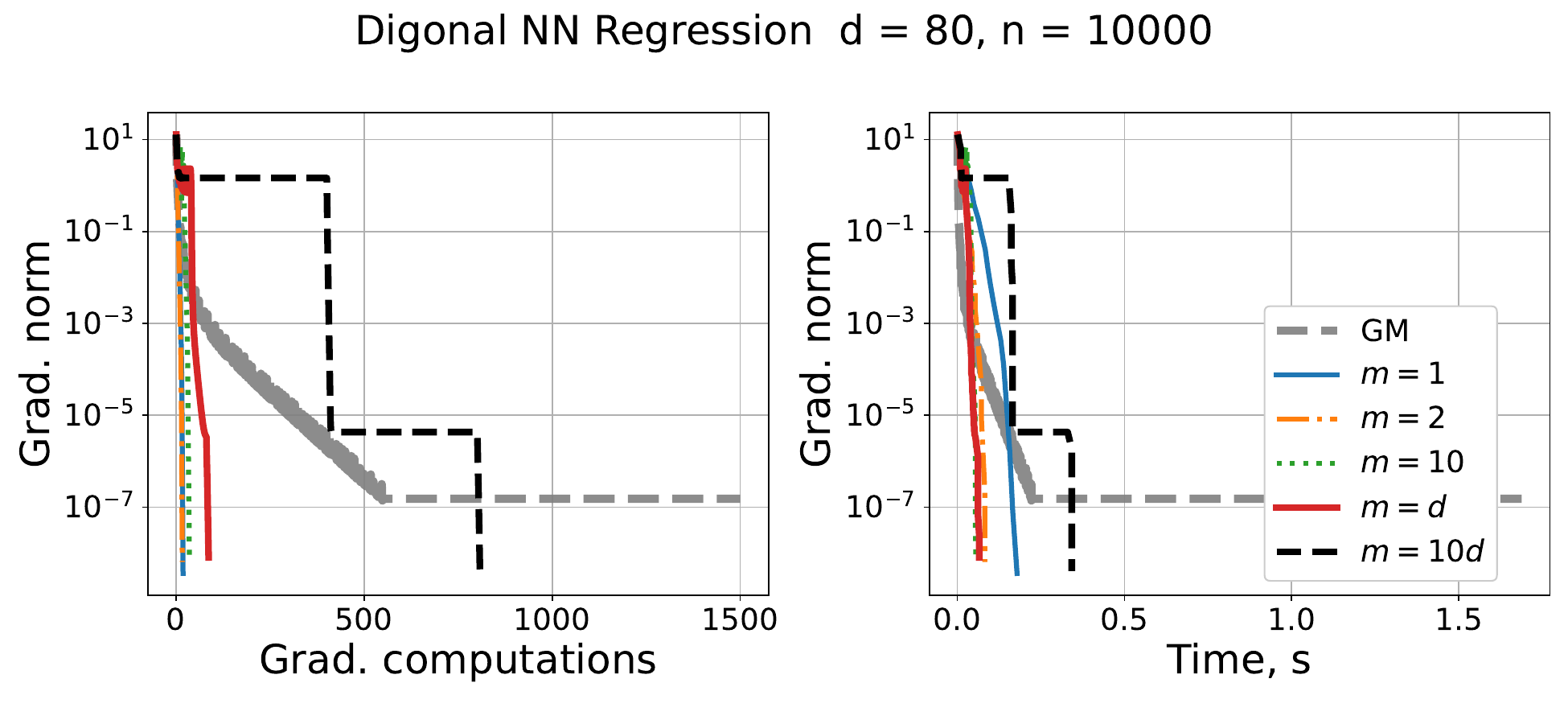} \\
	\includegraphics[width=0.42\textwidth ]{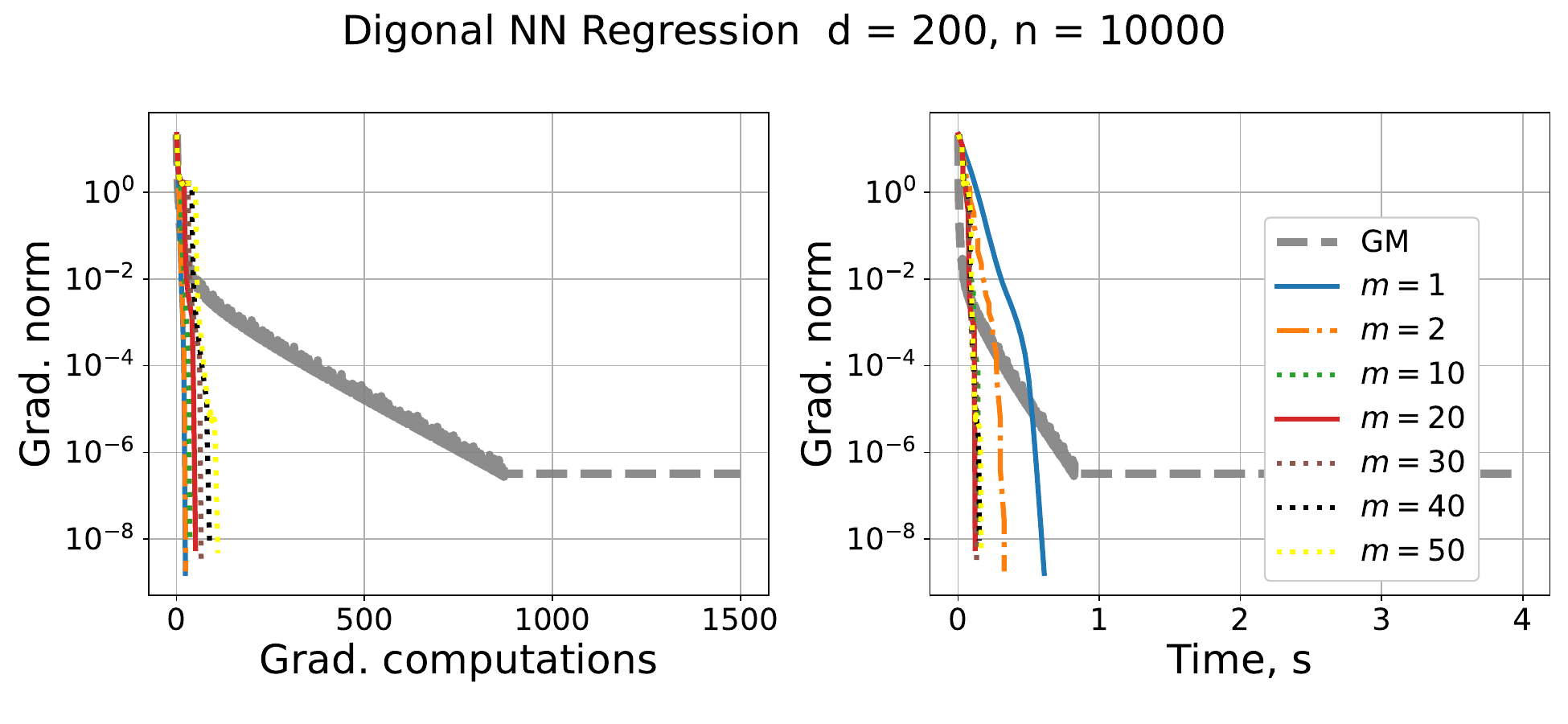} 
	\hspace*{-5pt}
	\caption{\small Performance of a Diagonal Neural Network trained on random data. }
	\label{fig:DNN}
\end{figure}

Finally, we consider non-convex optimization problem with the following objective:
$$
\ba{rcl}
f(\xx,\yy) & = & \| \mat{A} ( \xx \odot \yy) - \bb  \|^2, \qquad \xx, \yy \in \R^d, 
\ea
$$
where $\mat{A} \in \R^{n \times d}$, $\bb \in \R^n$ are generated randomly with
the entries from Gaussian distribution, and $\odot$ denotes the coordinate-wise product.
The results are shown in Figure~\ref{fig:DNN}. We see that there are many choices of parameter $m > 1$ (the frequency of the Hessian updates) that lead to a considerable time saving without hurting the convergence rate. Notice that $m = d$ is not necessarily optimal for this problem. Our intuition is that the nature of the parametrization of our Diagonal Neural Network implies a certain \textit{effective dimension} that is smaller than the full dimension $d$.
Hence, it can be an interesting direction for further research to develop
new strategies of the Hessian updates that are suitable for the problems with a specific structure.

\end{document}